\documentclass[12pt]{article}
\usepackage{latexsym,amssymb,amsfonts,amsmath,amsthm,graphicx}
\usepackage{setspace}

\def\beq{ \begin{equation} }
\def\eeq{ \end{equation} }
\def\beqx{ \begin{equation*} }
\def\eeqx{ \end{equation*} }
\def\beqa{\begin{eqnarray}}
\def\eeqa{\end{eqnarray}}
\def\beqax{\begin{eqnarray*}}
\def\eeqax{\end{eqnarray*}}
\def\mn{\medskip\noindent}
\def\ms{\medskip}

\def\square{\vcenter{\vbox{\hrule height .4pt
  \hbox{\vrule width .4pt height 5pt \kern 5pt
        \vrule width .4pt} \hrule height .4pt}}}
\def\eopt{\hfill$\square$}
\def\nass{\noalign{\smallskip}}
\def\sqz{\kern -0.15em}

\newcommand{\calN}{\mathcal{N}}
\newcommand{\ep}{\varepsilon}
\newcommand{\ape}{\simeq}
\numberwithin{equation}{section}
\newtheorem{theorem}{Theorem}
\newtheorem{conj}{Conjecture}
\newtheorem{lemma}{Lemma}
\newtheorem{remark}{Remark}

\begin{document}

\title{Evolutionary dynamics of tumor progression \\ with random fitness values}
\author{Rick Durrett$^{1,}$\thanks{Partially supported by NSF grant DMS 0704996 from the probability program.} , Jasmine Foo$^{2,}$\thanks{Partially supported by NIH grant R01CA138234.}, Kevin Leder$^{2,}$\thanks{Partially supported by NIH grant U54CA143798.}, \\ John Mayberry$^{1,}$\thanks{Corresponding Author. Email: jm858@cornell.edu, Tel.: +1 607 255 8262, Fax: +1 607 255 7149.} \thanks{Partially supported by NSF RTG grant DMS 0739164.} , and Franziska Michor$^{2,}$\thanks{Partially supported by NIH grants R01CA138234 and U54CA143798, a Leon Levy Foundation Young Investigator Award, and a Gerstner Young Investigator Award.}   \\ \\
$^1$ Department of Mathematics, Cornell University, Ithaca, NY 14853\\
$^2$ Computational Biology Program, Memorial Sloan-Kettering \\ Cancer Center, New York, NY 10065}
\maketitle

\begin{abstract}

Most human tumors result from the accumulation of multiple genetic and
epigenetic
alterations in a single cell. Mutations that confer a fitness advantage to
the cell are known as driver mutations and are causally related to
tumorigenesis. Other mutations, however, do not change the phenotype of the
cell or even decrease cellular fitness. While much experimental effort is
being devoted to the identification of the different functional effects of
individual mutations, mathematical modeling of tumor progression generally
considers constant fitness increments as mutations are accumulated. In this
paper we study a mathematical model of tumor progression with random
fitness increments. We analyze a multi-type branching process in which
cells accumulate mutations whose fitness effects are chosen from a
distribution. We determine the effect of the fitness distribution on the
growth kinetics of the tumor. This work contributes to a quantitative
understanding of the accumulation of mutations leading to cancer phenotypes.
\end{abstract}

\bigskip


\newpage

\section{Introduction}

Tumors result from an evolutionary process occurring within a tissue (Nowell, 1976). From an evolutionary point of view, tumors can be considered as collections of cells that accumulate genetic and epigenetic alterations. The phenotypic changes that these alterations confer to cells are subjected to the selection pressures within the tissue and lead to adaptations such as the evolution of more aggressive cell types, the emergence of resistance, induction of angiogenesis, evasion of the immune system, and colonization of distant organs with metastatic growth. Advantageous heritable alterations can cause a rapid expansion of the cell clone harboring such changes, since these cells are capable of outcompeting cells that have not evolved similar adaptations. The investigation of the dynamics of cell growth, the speed of accumulating mutations, and the distribution of different cell types at various timepoints during tumorigenesis is important for an understanding of the natural history of tumors. Further, such knowledge aids in the prognosis of newly diagnosed tumors, since the presence of cell clones with aggressive phenotypes lead to less optimistic predictions for tumor progression. Finally, a knowledge of the composition of tumors allows for the choice of optimum therapeutic interventions, as tumors harboring pre-existing resistant clones should be treated differently than drug-sensitive cell populations.

Mathematical models have led to many important insights into the dynamics of tumor progression and the evolution of resistance (Goldie and Coldman, 1983 and 1984; Bodmer and Tomlinson, 1995; Coldman and Murray, 2000; Knudson, 2001; Maley and Forrest, 2001; Michor et al., 2004; Iwasa et al., 2005; Komarova and Wodarz, 2005; Michor et al., 2006; Michor and Iwasa, 2006; Frank 2007; Wodarz and Komarova, 2007). These mathematical models generally fall into one of two classes: (i) constant population size models, and (ii) models describing exponentially growing populations. Many theoretical investigations of exponentially growing populations employ multi-type branching process models (e.g., Iwasa et al., 2006; Haeno et al., 2007; Durrett and Moseley, 2009), while others use population genetic models for homogeneously mixing exponentially growing populations (e.g., Beerenwinkel et al., 2007; Durrett and Mayberry, 2009). In this paper, we focus on branching process models. In these models, cells with $i\ge 0$ mutations are denoted as type-$i$ cells, and $Z_i(t)$ specifies the number of type-$i$ cells at time $t$. Type-$i$ cells die at rate $b_i$, give birth to one new type-$i$ cell at rate $a_i$, and give birth to one new type-$(i+1)$ cell at rate $u_{i+1}$. In an alternate version, mutations occur with probability $\mu_{i+1}$ during birth events which occur at rate $\alpha_i$. These two versions are equivalent provided $u_{i+1}=\alpha_i\mu_{i+1}$ and $a_i = \alpha_i(1-\mu_{i+1})$. However, the relationship between the parameters must be kept in mind when comparing results between different formulations of the model.

One biologically unrealistic aspect of this model as presented in the literature is that all type-$i$ cells are assumed to have the same birth and death rates. This assumption describes situations during tumorigenesis in which the order of mutations is predetermined, i.e. the genetic changes can only be accumulated in a particular sequence and all other combinations of mutations lead to lethality. Furthermore, in this interpretation of the model, there cannot be any variability in phenotype among cells with the same number of mutations. In many situations arising in biology, however, there is marked heterogeneity in phenotype even if genetically, the cells are identical (Elowitz et al., 2002; Becskei et al., 2005; Kaern et al., 2005; Feinerman et al., 2008). This variability may be driven by stochasticity in gene expression or in post-transcriptional or post-translational modifications. In this paper, we modify the branching process model so that mutations alter cell birth rates by a random amount.

An important consideration for this endeavor is the choice of the mutational fitness distribution. The exponential distribution has become the preferred candidate in theoretical studies of the genetics of adaptation. The first theoretical justification of this choice was given by Gillespie (1983, 1984), who argued that if the number of possible alleles is large and the current allele is close to the top of the rank ordering in fitness values, then extreme value theory should provide insight into the distribution of the fitness values of mutations. For many distributions including the normal, Gamma, and lognormal distributions, the maximum of $n$ independent draws, when properly scaled, converges to the Gumbel or double exponential distribution, $\Lambda(x) = \exp(-e^{-x})$. In the biological literature, it is generally noted that this class of distributions only excludes exotic distributions like the Cauchy distribution, which has no moments. However, in reality, it eliminates all distributions with $P(X>x) \sim Cx^{-\alpha}$. For distributions in the domain of attraction of the Gumbel distribution, and if $Y_1 > Y_2 \cdots > Y_k$ are the $k$ largest observations in a sample of size $n$, then there is a sequence of constants $b_n$ so that the spacings $Z_i = i(Y_i-Y_{i+1})/b_n$ converge to independent exponentials with mean 1, see e.g., Weissman (1978). Following up on Gillespie's work, Orr (2003) added the observation that in this setting, the distribution of the fitness increases due to beneficial mutations has the same distribution as $Z_1$ independent of the rank $i$ of the wild type cell.

To infer the distribution of fitness effects of newly emerged beneficial mutations, several experimental studies were performed; for examples, see Imhoff and Schlotterer (2001), Sanjuan et al. (2004), and Kassen and Bataillon (2006). The data from these experiments is generally consistent with an exponential distribution of fitness effects. However, there is an experimental caveat that cannot be neglected (Rozen et al., 2002): if only those mutations are considered that reach 100\% frequency in the population, then the exponential distribution is multiplied by the fixation probability. By this operation, a distribution with a mode at a positive value develops. In a study of a quasi-empirical model of RNA evolution in which fitness was based on secondary structures, Cowperthwaite et al. (2005) found that fitnesses of randomly selected genotypes appeared to follow a Gumbel-type distribution. They also discovered that the fitness distribution of beneficial mutations appeared exponential only when the vast majority of small-effect mutations were ignored. Furthermore, it was determined that the distribution of beneficial mutations depends on the fitness of the parental genotype (Cowperthwaite et al., 2005; MacLean and Buckling, 2009). However, since the exceptions to this conclusion arise when the fitness of the wild type cell is low, these findings do not contradict the picture based on extreme value theory.

In contrast to the evidence above, recent work of Rokyta et al. (2008) has shown that in two sets of beneficial mutations arising in the bacteriophage ID11 and in the phage $\phi 6$ -- for which the mutations were identified by sequencing -- beneficial fitness effects are not exponential. Using a statistical method developed by Biesal et al.~(2007), they tested the null hypothesis that the fitness distribution has an exponential tail. They found that the null hypothesis could be rejected in favor of a distribution with a right truncated tail. Their data also violated the common assumption that small-effect mutations greatly outnumber those of large effect, as they were consistent with a uniform distribution of beneficial effects. A possible explanation for the bounded fitness distribution may be found in the culture conditions utilized in the experiments: they evolved ID11 on {\it E.coli} at an elevated temperature ($37^o$ C instead of $33^o$ C). There may be a limited number of mutations that will enable ID11 to survive in increased temperatures. The latter situation may be similar to scenarios arising during tumorigenesis, where, in order to develop resistance to a drug or to progress to a more aggressive stage, the conformation of a particular protein must be changed or a certain regulatory network must be disrupted. If there is a finite, but large, number of possible beneficial mutations, then it is convenient to use a continuous distribution as an approximation.

In this paper, we consider both bounded distributions and unbounded distributions for the fitness advance and derive asymptotic results for the number of type-$k$ individuals at time $t$. We determine the effects of the fitness distribution on the growth kinetics of the population, and investigate the rates of expansion for both bounded and unbounded fitness distributions. This model provides a framework to investigate the accumulation of mutations with random fitness effects.

The remainder of this section is dedicated to statements and discussion of our main results. Proofs of these results can be found in Sections \ref{prelimsec}-\ref{unbndsec}.

\subsection{Bounded distributions} \label{bndintro}

Let us consider a multi-type branching process in which type-$i$ cells have accumulated $i \geq 0$ advantageous mutations. Suppose the initial population consists entirely of type-0 cells that give birth at rate $a_0$ to new type-0 cells, die at rate $b_0 < a_0$, and give birth to new type-1 cells at rate $u_1$. The parameters $a_0$, $b_0$, and $u_1$ denote the birth rate, death rate, and mutation rate for type-0 cells. To simplify computations, we will approximate the number of type-0 cells by $Z_0(t)=V_0e^{\lambda_0t}$, where $\lambda_0 = a_0 - b_0 >0$. If the initial cell population $Z_0(0) = V_0 \gg 1/\lambda_0$, then the branching process giving the number of 0's is almost deterministic and this approximation is accurate. When a new type-1 cell is born, we choose $x > 0$ according to a continuous probability distribution $\nu$. The new type 1-cell and its descendants then have birth rate $a_0 + x$, death rate $b_0$, and mutation rate $u_2$. In general, type-$k$ cells with birth rate $a$ mutate to type-$(k+1)$ cell at rate $u_{k+1}$ and when a mutation occurs, the new type-$(k+1)$ cell and its descendants have an increased birth rate $a+ x$ where $x >0$ is drawn according to $\nu$. We let $Z_k(t)$ denote the total number of type-$k$ cells in the population at time $t$. When we refer to the $k$th generation of mutants, we mean the set of all type-$k$ cells.

We begin by considering situations in which the distribution of the increase in the birth rate is concentrated on $[0,b]$. In particular, suppose that $\nu$ has density $g$ with support in $[0,b]$ and assume that $g$ satisfies:
$$
\hbox{$g$ is continuous at $b$, $g(b)>0$, $g(x) \le G$ for $x \in [0,b]$}
\leqno(\ast)
$$
Our first result describes the mean number of first generation mutants at time $t$, $E Z_1(t)$.

\begin{theorem} \label{th1A} If $(\ast)$ holds, then
$$
EZ_1(t) \sim \frac{V_0 u_1 \, g(b)}{bt} e^{(\lambda_0 + b) t}
$$
where $a(t) \sim b(t)$ means $a(t)/b(t) \to 1$.
\end{theorem}


The next result shows that the actual growth rate of type-1 cells is slower than the mean. Here, and in what follows, we use $\Rightarrow$ to indicate convergence in distribution.

\begin{theorem} \label{th1B}
If $(\ast)$ holds and $p = b/\lambda_0$, then for $\theta\ge 0$,
\begin{align} \label{LTlim}
E\exp( - \theta t^{1+p} e^{-(\lambda_0+b)t} Z_1(t) ) & \to \exp( - V_0u_1 \theta^{\lambda_0/(\lambda_0+b)} c_1(\lambda_0,b) ),
\end{align}
where $c_1(\lambda_0,b)$ is an explicit constant whose value will be given in \eqref{c1def}. In particular, we have
$$
t^{1+p} e^{-(\lambda_0+b)t} Z_1(t) \Rightarrow V_1,
$$
where $V_1$ has Laplace transform given by the righthand side of \eqref{LTlim}.
\end{theorem}

\mn
Theorem \ref{th1B} is similar to Theorem 3 in Durrett and Moseley (2009) which assumes a deterministic fitness distribution so that all type-1 cells have growth rate $\lambda_1=\lambda_0+b$. There, the asymptotic growth rate of the first generation is $\exp(\lambda_1 t)$. In contrast, the continuous fitness distribution we consider here has the effect of slowing down the growth rate of the first generation by the polynomial factor $t^{1+p}$. To explain this difference, we note that the calculation of the mean given in Section \ref{bndZ1sec} shows that the dominant contribution to $Z_1(t)$ comes from growth rates $x=b-O(1/t)$. However, mutations with this growth rate are unlikely until the number of type-0 cells is $O(t)$, i.e., roughly at time $r_1 = (1/\lambda_0)\log t$. Thus at time $t$, the number of type-1 cells will be roughly $\exp((\lambda_0+b)(t-r_1)) = \exp((\lambda_0+b)t)/t^{1+p}$.

\begin{figure}[tbp] 
  \centering
  \includegraphics[width=4.25in,height=3.5in,keepaspectratio]{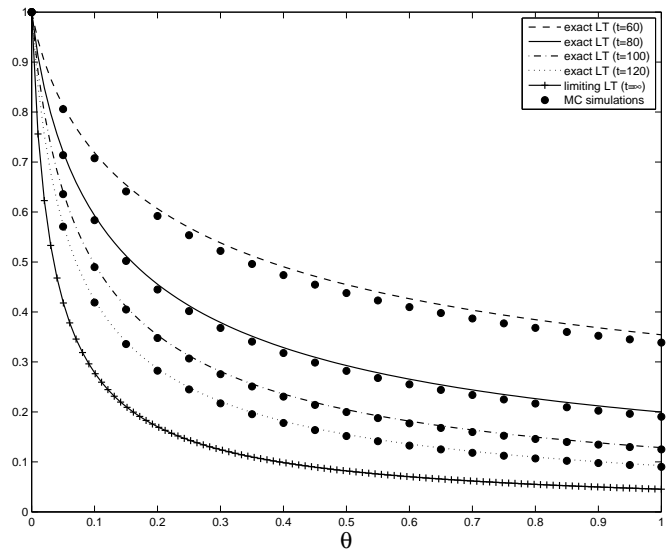}
  \caption{Plot of the exact Laplace transform (LT) for $t^{(1+p)}e^{-(\lambda_0+b)t} Z_1(t)$ at times $t=60,80,100,120$, the approximations from Monte Carlo (MC) simulations at the corresponding times, and the asymptotic Laplace transform from Theorem \ref{th1B}. Parameter values: $a_0 = 0.2$, $b_0=0.1$, $b = 0.01$, and $u_1=10^{-3}$. $g$ is uniform on $[0,.01]$.}
  \label{fig:wave1}
\end{figure}

To prove Theorem \ref{th1B}, we look at mutations as a point process in $[0,t] \times [0,b]$: there is a point at $(s,x)$ if there was a mutant with birth rate $a_0+x$ at time $s$. This allows us to derive the following explicit expression for the Laplace transform of $Z_1(t)$:
\beqx
E( e^{-\theta Z_1(t)} ) = \exp\left( - u_1 \int_{0}^b dx \, g(x) \int_0^t ds \, V_0 e^{\lambda_0 s} ( 1 - \tilde\phi_{x,t-s}(\theta) ) \right)
\eeqx
where $\tilde\phi_{x,r}(\theta) = Ee^{-\theta \tilde Z^x_r}$ and $\tilde Z^x_r$ is a continuous-time branching process with birth rate $a_0+x$, death rate $b_0$, and initial population $\tilde{Z}_0^x = 1$. In Figure \ref{fig:wave1}, we compare the exact Laplace transform of $t^{1+p} \exp(-(\lambda_0+b)t) Z_1(t)$ with the results of simulations and the limiting Laplace transform from Theorem \ref{th1B}, illustrating the convergence as $t \to \infty$.

Notice that the Laplace transform of $V_1$ has the form $\exp(C \, \theta^\alpha)$ where $\alpha = \lambda_0/(\lambda_0+b)$ which implies that $P(V_1 > v) \sim v^{-\alpha}$ as $v \to \infty$ (see, for example, the argument in Section 3 of Durrett and Moseley (2009)). To gain some insight into how this limit comes about, we give a second proof of the convergence that tells us the limit is the sum of points in a nonhomogeneous Poisson process. Each point in the limiting process represents the contribution of a different mutant lineage to $Z_1(t)$.

\begin{theorem} \label{th1C}
$V_1 = \lim_{t\to\infty} t^{1+p} e^{-(\lambda_0+b)t} Z_1(t)$ is the sum of the points of a Poisson process
on $(0,\infty)$ with mean measure $\mu(z,\infty) = A_1(\lambda_0,b) u_1 V_0 z^{-\lambda_0/(\lambda_0+b)}$.
\end{theorem}

\noindent
A similar result can be obtained for deterministic fitness distributions, see the Corollary to Theorem 3 in Durrett and Moseley (2009). However, the new result shows that the point process limit is not an artifact of assuming that all first generation mutants have the same growth rate. Even when the fitness advances are random, different mutant lines contribute to the limit. This result is consistent with observations of Maley et al.~(2006) and Shah et al.~(2009) that tumors contain cells with different mutational haplotypes. Theorem \ref{th1C} also gives quantitative predictions about the relative contribution of different mutations to the total population. These implications will be explored further in a follow-up paper currently in progress.

With the behavior of the first generation analyzed, we are ready to proceed to the study of further generations. The computation of the mean is straightforward.

\begin{theorem} \label{th2A}
If $(\ast)$ holds, then
$$
EZ_k(t) \sim \frac{V_0 \cdot u_1 \cdots u_k \cdot g(b)^k} {t^k b^k k!} e^{(\lambda_0 + kb) t}
$$
\end{theorem}

\noindent
As in the $k=1$ case, the mean involves a polynomial correction to the exponential growth and again, does not give the correct growth rate for the number of type-$k$ cells. To state the correct limit theorem describing the growth rate of $Z_k(t)$, we will define $p_k$ and $u_{1,k}$ by
$$
k+p_k = \sum_{j=0}^{k-1} \frac{\lambda_0+kb}{\lambda_0+jb} \quad\hbox{and}\quad
u_{1,k}  = \prod_{j=1}^k u_j^{\lambda_0/(\lambda_0+(j-1)b)}
$$
for all $k \geq 1$.

\begin{theorem} \label{th2B}
If $(\ast)$ holds, then for $\theta\ge 0$
\begin{align*}
E\exp( - \theta t^{k+p_k} e^{-(\lambda_0+kb)t} Z_k(t) ) & \to \exp( - c_k(\lambda_0,b)
V_0u_{1,k} \theta^{\lambda_0/(\lambda_0+kb)} ) \\
t^{k+p_k} e^{-(\lambda_0+kb)t} Z_k(t) & \Rightarrow V_k
\end{align*}
\end{theorem}

We prove this result by looking at the mutations to type-1 individuals as a three dimensional Poisson point process: there is a point at $(s,x,v)$ if there was a type-1 mutant with birth rate $a_0+x$ at time $s$ and the number of its type-1 descendants at time
$t$, $Z_1^{s,x}(t)$, has $e^{-(\lambda_0+x)(t-s)} Z_1^{s,x}(t) \to v$ with $v>0$. To study $Z_k(t)$ we will
let $Z_k^{s,x,v}(t)$ be the type-$k$ descendants at time $t$ of the 1 mutant at $(s,x,v)$. $Z_k^{s,x,v}$
is the same as a process in which the initial type (here type-1 cells) behaves like
$ve^{(\lambda_0+x)(t-s)}$ instead of $Z_0(t) = V_0 e^{\lambda_0t}$, so the result can be proved by induction.

To explain the form of the result we consider the case $k=2$. Breaking things down according to the times and the sizes of the mutational changes, we have
\begin{align*}
EZ_2(t) = \int_0^b dx_1 \, g(x_1) \int_0^b dx_2\, g(x_2) & \int_0^t ds_1 \int_{s_1}^t ds_2 \\
& V_0 e^{\lambda_0 s_1} u_1 e^{(\lambda_0+x_1) (s_2-s_1)} u_2 e^{(\lambda_0+x_1+x_2) (t-s_2)}
\end{align*}
As in the result for $Z_1(t)$ the dominant contribution comes from $x_1, x_2 = b - O(1/t)$ and as in the discussion preceding the statement of Theorem \ref{th1B}, the time of the first mutation to $b-O(1/t)$ is $ \approx r_1 = (\log t)/\lambda_0$. The descendants of this mutation grow at exponential rate $\lambda_0+b-O(1/t)$, so the time of the first mutation to $2b-O(1/t)$ is $\approx r_2 = r_1 + (\log t)/(\lambda_0+b)$. Noticing that
$$
\exp( (\lambda_0+2b) (t-r_1-r_2)) = \exp((\lambda_0+2b)t) t^{-(\lambda_0+2b)/\lambda_0 - (\lambda_0+2b)/(\lambda_0+b)}
$$
tells us what to guess for the polynomial term: $t^{-(2+p_2)}$ where
$$
2 + p_2 = \frac{\lambda_0 + 2b}{\lambda_0}  + \frac{\lambda_0 + 2b}{\lambda_0+b}
$$

\begin{figure}[tbp] 
  \centering
  \includegraphics[width=4.25in,height=3.5in,keepaspectratio]{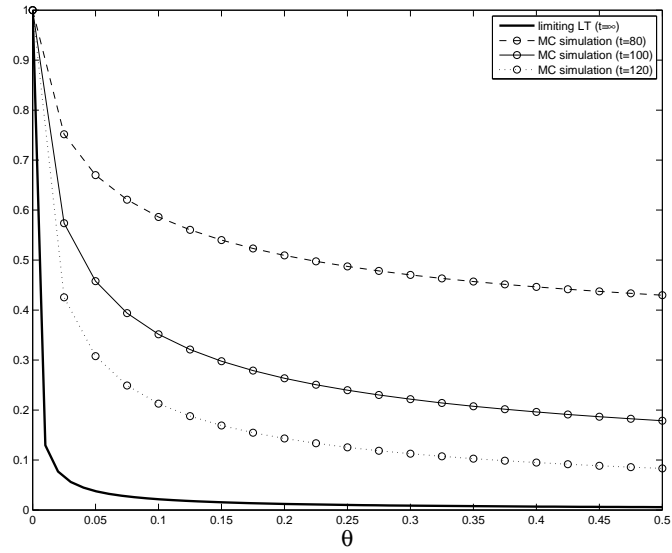}
  \caption{Plot of the approximations to the Laplace transform of $t^{2+p_2} e^{-(\lambda_0+2b)t} Z_2(t)$ from Monte Carlo (MC) simulations at times $t=80,100,120$ along with the asymptotic Laplace transform from Theorem \ref{th2B}. Parameter values: $a_0 = 0.2$, $b_0=0.1$, $b = 0.01$, and $u_1=u_2 = 10^{-3}$. $g$ is uniform on [0,0.01].}
  \label{fig:wave2}
\end{figure}

In Figure \ref{fig:wave2}, we compare the asymptotic Laplace transform from Theorem \ref{th2B} with the results of simulations in the case $k=2$. To explain the slow convergence to the limit, we note that if we take account of the mutation rates $u_1,u_2$ in the heuristic from the previous paragraph (which becomes important when $u_1, u_2$ are small), then the first time we see a type-1 cell with growth rate $b-O(1/t)$ will not occur until time $\lambda_0^{-1} \log (t/u_1)$ when the type-0 cells reach $O(t/u_1)$ and so the first type-2 cell with growth rate $2b - O(1/t)$ will not be born until time $r = \lambda_0^{-1} \log (t/u_1) + (\lambda_0+b)^{-1} \log (t/u_2)$ when the descendants of the type-1 cells with growth rate $b - O(1/t)$ reach size $O(t/u_2)$. When $u_1=u_2=10^{-3}$, $\lambda_0 = .1$, and $b=.01$, $r \approx 223$. The mutations created at this point will need some time to grow and become dominant in the population. It would be interesting to compare simulations at time 300, but we have not been able to do this due to the large number of different growth rates in generation 1.

\subsection{Unbounded distributions} \label{unbndintro}

Let us now consider situations in which the fitness distribution is unbounded. Suppose that the fitness increase follows a generalized Frechet distribution,
\beq \label{genunbnd}
P(X > x) = x^\beta e^{-\gamma x^\alpha}
\eeq
for some positive $\gamma, \alpha$ and any $\beta \in \mathbb{R}$. There is a two-fold purpose for considering such distributions. First, if i.i.d.~random variables $\zeta_1, \ldots, \zeta_n$ have a power law tail, i.e. $P(\zeta_i > y) \sim cy^{-\alpha}$ as $y \to \infty$, then their maxima and the spacings between order statistics converge to a limit of the form \eqref{genunbnd} with $\beta=0$. Second, this choice allows us to consider the gamma($\beta+1,\gamma$) distribution which has $\alpha=1$ and the normal distribution, which asymptotically has this form with $\alpha= 2, \beta=-1$.

To analyze this situation, we will again take a Poisson process viewpoint and look at the contribution from a mutation at time $s$ with increased growth rate $x$. A mutation that increases the growth rate by $x$ at time $s$ will, if it does not die out, grow to $e^{(\lambda_0+x)(t-s)} \zeta$ at time $t$ where $\zeta$ has an exponential distribution. The growth rate $(\lambda_0+x)(t-s) \ge z$ when
$$
x \ge \frac{z}{t-s} - \lambda_0.
$$
Therefore,
\begin{align*}
\mu(z,\infty) & \equiv E( \hbox{\# mutations with }(\lambda_0+x)(t-s) \ge z ) \\
& = V_0 u_1 \int_0^t \left( \frac{z}{t-s} - \lambda_0 \right)^\beta e^{\lambda_0s}
\exp\left( - \gamma \left( \frac{z}{t-s} - \lambda_0 \right)^{\sqz\alpha} \, \right) \, ds \\
&= V_0 u_1 \int_0^t \left( \frac{z}{t-s} - \lambda_0 \right)^\beta \exp(\phi(s,z)) \, ds
\end{align*}
where
\beq \label{maxme}
\phi(s,z) = \lambda_0 s - \gamma \left( \frac{z}{t-s} - \lambda_0 \right)^{\sqz\alpha}.
\eeq
The size of this integral can be found by maximizing the exponent $\phi$ over $s$ for fixed $z$. Since
\beq
\frac{\partial \phi}{\partial s}(s,z)
=  \lambda_0  -  \alpha \gamma \left( \frac{z}{t-s} - \lambda_0 \right)^{\sqz\alpha-1} \frac{z}{(t-s)^2}
\label{der1}
\eeq
and
\beq
\frac{\partial^2 \phi}{\partial s^2}(s,z)
= -  \alpha(\alpha-1) \gamma \left( \frac{z}{t-s} - \lambda_0 \right)^{\sqz\alpha-2} \frac{z^2}{(t-s)^4}
-  \alpha \gamma \left( \frac{z}{t-s} - \lambda_0 \right)^{\sqz\alpha-1} \frac{2z}{(t-s)^3}
\label{der2}
\eeq
we can see that $\partial^2 \phi/\partial s^2(s,z) < 0$ when $\alpha z> \lambda_0 (t-s) $ so that for all $z$ in this range, $\phi(s,z)$ is concave as a function of $s$ and achieves its maximum at a unique value $s_z$.

When $\alpha=1$, it is easy to set \eqref{der1} to 0 and solve for $s_z$. This in turn leads to an asymptotic formula for $\mu(z,\infty)$ and allows us to derive the following limit theorem for $Z_1(t)$.

\begin{theorem} \label{uthexp}
Suppose $\alpha =1$ and let $c_0 = \lambda_0/4\gamma$. Then $t^{-2} \log Z_1(t) \to c_0$ and
\beqx
\frac{1}{t} \left[ \log Z_1(t) - c_0 t^2 \left(1  +  \frac{(2\beta+1)\log t}{\lambda_0t} \right) \right] \Rightarrow y^*
\eeqx
where $y^*$ is the rightmost point in the point process with intensity given by
\beq \label{ptp}
(2c_0)^\beta  ( \pi/\lambda_0)^{1/2} V_0 u_1 \exp( \gamma\lambda_0 - \lambda_0 y/2c_0 ).
\eeq
\end{theorem}

When $\alpha \neq 1$, solving for $s_z$ becomes more difficult, but we are still able to prove the following limit theorem for $Z_1(t)$.

\begin{theorem} \label{uthW}
Suppose $\alpha>1$ is an integer.
There exist explicitly calculable constants $c_k = c_k (\alpha,\gamma)$, $0 \leq k < \alpha$, and $\kappa = \kappa(\beta,\alpha,\gamma)$ so that $t^{-(\alpha+1)/\alpha} \log Z_1(t) \to c_0$ and
\beqx
\frac{1}{t^{1/\alpha}} \left[ \log Z_1(t) - c_0 t^{(\alpha+1)/\alpha}\left(1 +  \sum_{1 \leq k < \alpha } c_k t^{-k/\alpha} + \kappa \frac{\log t}{t}\right)\right] \Rightarrow y^*
\eeqx
where $y^*$ is the rightmost particle in a point process with explicitly calculable intensity.
\end{theorem}

\noindent
The complicated form of the result is due to the fact that the fluctuations are only of order $t^{1/\alpha}$, so we have to be very precise in locating the maximum. The explicit formulas for the constants and the intensity of the point process are given in \eqref{cis} and \eqref{intensity2}. With more work this result could be proved for a general $\alpha>1$, but we have not tried to do this or prove Conjecture 1 below because the super-exponential growth rates in the unbounded case are too fast to be realistic.


We conclude this section with two comments. First, the proof of Theorem \ref{uthW} shows that in contrast to the bounded case, in the unbounded case, most type-1 individuals are descendants of a single mutant. Second, the proof shows that the distribution of the mutant with the largest growth rate is born at time $s \sim t /(\alpha+1)$ (see Remark \ref{timeofb} at the end of Section \ref{unbndsec}) and has growth rate $z = O(t^{(\alpha+1)/\alpha})$. The intuition behind this is that since the type-0 cells have growth rate $e^{\lambda_0 s}$ and the distribution of the increase in fitness has tail $\approx e^{-\gamma x^{\alpha}}$, the largest advance $x$ attained by time $t$ should occur when $s = O(t)$ and satisfy
$$
e^{C \lambda_0 t} e^{-\gamma x^\alpha}  =  O(1) \quad \text{ or } \quad  x = O(t^{1/\alpha}).
$$
The growth rate of its family is then $(\lambda_0 + x)(t-s) = O(t^{(\alpha+1)/\alpha})$.

Since the type-1 cells grow at exponential rate $c_1 t^{(\alpha+1)/\alpha}$,  if we apply this same reasoning to type-2 mutants, then the largest additional fitness advance $x$ attained by type-2 individuals should satisfy
$$
e^{c_1 t^{(\alpha+1)/\alpha}} e^{-\gamma x^\alpha}  =  O(1) \quad \text{ or } \quad x = O(t^{1/\alpha + 1/\alpha^2}).
$$
and the growth rate of its family will be $O(t^{1 + 1/\alpha + 1/\alpha^2})$. Extrapolating from the first two
generations, we make the following

\begin{conj}
Let $q(k) = \sum_{j=0}^k \alpha^{-j}$. As $t\to\infty$,
$$
\frac{1}{t^{q(k)}} \log Z_k(t) \to c_k
$$
\end{conj}

\noindent
Note that in the case of the exponential distribution, $q(k)=k+1$.

The rest of the paper is organized as follows. Sections \ref{prelimsec}-\ref{unbndsec} are devoted to proofs of our main results. After some preliminary notation and definitions in Section \ref{prelimsec}, Theorems \ref{th1A}-\ref{th1C} are proved in Section \ref{bndZ1sec}, Theorems \ref{th2A}-\ref{th2B} in Section \ref{bndZksec}, and Theorems \ref{uthexp}-\ref{uthW} in Section \ref{unbndsec}. We conclude with a discussion of our results in Section \ref{discsec}.

\section{Preliminaries} \label{prelimsec}
This section contains some preliminary notation and definitions which we will need for the proofs of our main results. We denote by $\calN(t)$ the points in a two dimensional Poisson process on $[0,t] \times [0,\infty)$ with mean measure
$$
V_0e^{\lambda_0 s}  ds \nu(dx)
$$
where in Sections \ref{bndZ1sec}-\ref{bndZksec}, $\nu(dx) = g(x)dx$ with $g$ satisfying $(\ast)$ and in Section \ref{unbndsec}, $\nu$ has tail $\nu(x,\infty) =  x^{\beta} e^{-\gamma x^\alpha}$. In other words, we have a point at $(s,x)$ if there was a mutant with birth rate $a_0+x$ at time $s$. Define a collection of independent birth/death branching processes $Z_1^{s,x}(t)$ indexed by $(s,x) \in \calN(t)$ with $Z_1^{s,x}(s) = 1$, individual birth rate $a_0 +x$, and death rate $b$. $Z_1^{s,x}(t)$ is the contribution of the mutation at $(s,x)$ and
\beqx
Z_1(t) = \sum_{(s,x) \in \calN(t)} Z_1^{s,x}(t).
\eeqx
It is well known that
$$
e^{-(\lambda_0+x)(t-s)}Z_1^{s,x}(t) \to \frac{b}{a_0 +x} \delta_0 + \frac{\lambda_0 + x}{a_0 + x} \zeta
$$
where $\zeta \sim \exp((\lambda_0 + x)/(a_0 +x))$ (see, for example, equation (1) in Durrett and Moseley (2009)). In several results, we shall make use of the three dimensional Poisson process $\mathcal{M}(t)$ on $[0,t] \times [0,\infty) \times (0,\infty)$ with intensity
$$
 V_0 e^{\lambda_0s} \nu(dx) \left(\frac{\lambda_0 + x}{a_0 + x} \right)^2 e^{-v (\lambda_0+x)/(a_0+x)} dv.
$$
In words, $(s,x,v) \in \mathcal{M}(t)$ if there was a mutant with birth rate $a_0 + x$ at time $s$ and the number of its descendants at time $t$, $Z_1^{s,x}(t)$, has $ Z_1^{s,x}(t) \sim v e^{(\lambda_0+x)(t-s)}$. It is also convenient to define the mapping $z: [0,\infty) \times [0,t] \to [0,\infty)$ which maps a point $(s,x) \in \calN(t)$ to the growth rate of the induced branching process if it survives: $z(s,x) = (\lambda_0+x)(t-s)$ and let
$$
\mu(A) = E |\{(s,x) \in \calN(t): z(s,x) \in A\}|
$$
for $A \subset [0,\infty)$.

We shall use $C$ do denote a generic constant whose value may change from line to line. We write $f(t) \sim g(t)$ if $f(t)/g(t) \to 1$ as $t \to \infty$ and $f(t) = o(g(t))$ is $f(t)/g(t) \to 0$. $f(t) \gg (\ll) g(t)$ means that $f(t)/g(t) \to \infty \, ( \text{ resp. } 0)$ as $t \to \infty$ and $f(t) = O(g(t))$ means $|f(t)| \leq C g(t)$ for all $t > 0$. We also shall use the notation $f(t) \ape g(t)$ if $f(t) = g(t) + o(1)$ as $t \to \infty$.

\section{Bounded distributions, ${\bf Z}_{\bf 1}$} \label{bndZ1sec}

In this section, we prove Theorems \ref{th1A} - \ref{th1C}.

\mn
{\it Proof of Theorem \ref{th1A}.} Mutations to 1's occur at rate $V_0 e^{\lambda_0s}$ so
\begin{align}
EZ_1(t) & = u_1 \int_0^t \int_0^b e^{(t-s)(\lambda_0+x)} g(x) \, dx V_0 e^{\lambda_0 s} \, ds
\nonumber\\
& = u_1 V_0 e^{\lambda_0 t} \int_0^b dx \, g(x) \int_0^t e^{(t-s) x} \, ds
\label{EZ1}\\
& =  u_1 V_0 e^{\lambda_0 t} \int_0^b dx \, g(x) \frac{ e^{tx} - 1}{x}.
\nonumber
\end{align}
We begin by showing that the contribution from $x \in [0, b - (1+k)\log t)/t]$ can be ignored for any $k \in [0,\infty)$. The Mean Value theorem implies that
\beq
\frac{e^{tx}-1}{x} \leq  t e^{tx}
\label{expineq}
\eeq
Using this and the fact that $\int_c^d t e^{tx} dx  \leq e^{td}$ for any $c < d$, we can see that
\beq
t^k e^{-bt}  \int_0^{ b-(1+k)(\log t)/t } dx \, g(x) \frac{ e^{tx} - 1 }{x}
\le G t^{k}e^{-(1+k)\log t} \to 0
\label{xtrunc}
\eeq
To handle the other piece of the integral, we take $k=1$ and note that
$$
\int_{b-(2\log t)/t }^b dx \, g(x)  \frac{ e^{tx} - 1}{x}
\sim \frac{g(b)}{b} e^{bt} \int_{b-2 \log t/t}^b e^{t(x-b)} \, dx
$$
After changing variables $y = (b-x)t$, $dx = -dy/t$, the last integral
$$
= \frac{1}{t} \int_0^{2\log t} e^{-y} \, dy \sim 1/t
$$
which proves the result. \eopt

The above proof tells us that the dominant contribution to the 1's come from mutations with fitness increase $x \geq b_t = b-2 \log t/t$. To describe the times at which the dominant contributions occur, let $S(t)= (2/b) \log \log t$. Then the contribution to the mean from $x\in[b_t,b]$ and $s \ge S(t) $ is by (\ref{EZ1})
\begin{align*}
& \le Gu_1 V_0 e^{(\lambda_0+b)t} \frac{2(\log t)}{t} \int_{S(t)}^\infty e^{-sb_t} \, ds \\
& \le Gu_1 V_0 e^{(\lambda_0+b)t} \frac{2(\log t)}{tb_t} e^{-b_t S(t)}
\end{align*}
Since $b_tS(t) \geq 2 \log \log t$, this quantity is $o(t^{-1} e^{(\lambda_0+b)t})$. In words, the dominant contribution to the mean comes from points close to $(0,b)$ or more precisely from $[0, (2/b) \log \log t] \times [b - (2\log t)/t, b]$.

\begin{proof} [Proof of Theorem \ref{th1B}.] It suffices to prove \eqref{LTlim}. The computation in (\ref{xtrunc}) with $k=2+p$ implies that the contribution from mutations with $x \le b_t = b - (3+p)(\log t)/t$ can be ignored. Therefore, we have
$$
E \exp(-\theta Z_1(t)) \ape E \left( \exp(-\theta Z_1(t)); A_t \right)
$$
where $A_t = \{(s,x) \in \calN(t): x > b_t\}$. By Lemma 2 of Durrett and Moseley (2009), we have
$$
E( e^{-\theta Z_1(t)}; A_t ) = \exp\left( - u_1 \int_{b_t}^b dx \, g(x) \int_0^t ds \, V_0 e^{\lambda_0 s} ( 1 - \tilde\phi_{x,t-s}(\theta) ) \right)
$$
where $\tilde\phi_{x,r}(\theta) = Ee^{-\theta \tilde Z^x_r}$ and $\tilde Z^x_r$ is a birth/death branching process with birth rate $a_0+x$, death rate $b_0$, and initial population $\tilde{Z}_0^x = 1$. Using
\beq
e^{-(\lambda_0+b)t} = e^{-(\lambda_0+x)(t-s)} e^{-(\lambda_0+x)s} e^{-(b-x)t}
\label{expident}
\eeq
we have
\begin{align*}
E \left( \exp(-\theta Z_1(t) e^{-t(\lambda_0+b)} t^{1+p} ); A_t \right)
= \exp\biggl( - u_1 V_0 \int_{b_t}^b dx \, g(x)
\int_0^{t} ds \, e^{\lambda_0 s} & \\
 \qquad \{ 1 - \tilde\phi_{x,t-s}
(\theta e^{-(\lambda_0+x)(t-s)} e^{-(\lambda_0+x)s} e^{-(b-x)t} t^{1+p} ) \} \biggr) &
\end{align*}
Changing variables $s=r_x+r$ where $r_x = \frac{1}{\lambda_0+x}\log(t^{1+p})$ on the inside integral, $y=(b-x)t$, $dy/t = -dx$ on the outside, and continuing to write $x$ as short hand for $b-y/t$, the above
\begin{align} \label{exponen}
=  \exp\biggl( & - u_1 V_0 \int_0^{(3+p) \log t} \frac{dy}{t}
g(x) t^{(1+p) \lambda_0/(\lambda_0+x)} \nonumber \\
& \int_{-r_x}^{t-r_x} dr \, e^{\lambda_0 r} \{ 1 - \tilde\phi_{x,t-r-r_x}
(\theta e^{-(\lambda_0+x)(t-r-r_x)} e^{-(\lambda_0+x)r} e^{-y}  ) \} \biggr)
\end{align}
Formula (20) in Durrett and Moseley (2009) implies that as $u\to\infty$,
\beq
1 - \tilde\phi_{x,u}(\theta e^{-(\lambda_0+x)u}) \to \frac{\lambda_0+x}{a_0+x} \cdot
\frac{ \theta }{ \theta + \frac{\lambda_0+x}{a_0+x} }
\label{tildelim}
\eeq
and therefore, letting $t\to\infty$ and using $(1+p)\lambda_0/(\lambda_0+b)=1$, we can see that the expression in \eqref{exponen}
$$
\to \exp\left( - u_1 V_0 g(b) \int_0^{\infty} dy \,
\frac{\lambda_0+b}{a_0+b} \int_{-\infty}^{\infty} dr \, e^{\lambda_0 r}
\frac{ \theta e^{-(\lambda_0+b)r} e^{-y}}{ \theta e^{-(\lambda_0+b)r} e^{-y} + \frac{\lambda_0+b}{a_0+b} } \right)
$$
Changing variables $r=\frac{1}{\lambda_0+b}\{ q +  \log[\theta e^{-y} (a_0+b)/(\lambda_0+b)] \}$,
$dr = dq/(\lambda_0+b)$ gives
\begin{align*}
= \exp\biggl( - u_1 V_0 g(b) \theta^{\lambda_0/(\lambda_0+b)}
\left(\frac{\lambda_0+b}{a_0+b}\right)^{b/(\lambda_0+b)} & \int_0^{\infty} dy \, e^{-y\lambda_0/(\lambda_0+b)}\\
& \int_{-\infty}^{\infty} \frac{dq}{\lambda_0+b} \, e^{q\lambda_0/(\lambda_0+b)}
\frac{ e^{-q} }{ e^{-q} + 1 } \biggr)
\end{align*}
To simplify the first integral we note that
$$
\int_0^{\infty} dy \, e^{-y\lambda_0/(\lambda_0+b)} = \frac{\lambda_0+b}{\lambda_0}
$$
For the second integral, we prove

\begin{lemma} If $0<c<1$
\beq
\int_{-\infty}^{\infty} dq \, e^{qc}
\frac{ e^{-q} }{ e^{-q} + 1 }  = \Gamma(c)\Gamma(1-c)
\eeq
\end{lemma}

\begin{proof} We can rewrite the integral as
$$
\int_{-\infty}^\infty dq \, e^{qc} \int_0^\infty dx \, e^{-x} e^{-q} \exp(-e^{-q} x)
$$
so that after interchanging the order of integration and changing variables $w=e^{-q}x$, $dw = - dq \, e^{-q} x$ so that $w/x = e^{-q}$, $dw/x = - dq \, e^{-q}$, we have
$$
 =  \int_0^\infty dx \, \int_{0}^\infty \frac{dw}{x} (w/x)^{-c} e^{-x} e^{-w} \\
 = \int_0^\infty dx \, x^{-1+c} e^{-x}
\int_{0}^\infty dw \, w^{-c} e^{-w}
$$
which is $= \Gamma(c) \Gamma(1-c)$.
\end{proof}

Taking $c = \lambda_0/(\lambda_0+b)$ and letting
\beq
c_1(\lambda_0,b) = g(b) \frac{\lambda_0+b}{\lambda_0} \cdot \frac{1}{\lambda_0+b}
\left( \frac{ a_0+b }{\lambda_0+b} \right)^{-b/(\lambda_0+b)}
\Gamma(\lambda_0/(\lambda_0+b)) \Gamma(1-\lambda_0/(\lambda_0+b))
\label{c1def}
\eeq
we have proved Theorem \ref{th1B}. \end{proof}

Recall that we have assumed $Z_0(t) = V_0 e^{\lambda_0 t}$ is deterministic. This assumption can be relaxed to obtain the following generalization of Theorem \ref{th1B} which is used in Section \ref{bndZksec}.

\begin{lemma} \label{th1Bgen}
Suppose that $Z_0(t)$ is a stochastic process with $Z_0(t) \sim e^{\lambda_0 t} V_0$ for some constant $V_0$ as $t \to \infty$. Then the conclusions of Theorem \ref{th1B} remain valid.
\end{lemma}

To see why this is true, we can use a variant of Lemma 2 from Durrett and Moseley (2009) to conclude that
$$
E\left(e^{-\theta Z_1(t)}|\mathcal{F}_t^0\right)=\exp\left(-u_1\int_0^b dx \, g(x)\int_0^t ds Z_0(s)\left(1-\tilde{\phi}_{x,t-s}(\theta)\right)\right),
$$
where $\mathcal{F}_t^0$ is the $\sigma$-field generated by $Z_0(s)$ for $s\leq t$. Therefore,
$$
E\left(e^{-\theta Z_1(t)}\right)= E \exp\left(-u_1\int_0^b dx \, g(x)\int_0^t ds Z_0(s)\left(1-\tilde{\phi}_{x,t-s}(\theta)\right)\right),
$$
Given $\ep >0$, we can choose $t_\ep >0$ so that
$$\left|\frac{Z_0(t)}{V_0 \exp(\lambda_0 t)} - 1\right| < \ep$$
for all $t > t_\ep$. Since the contribution from $t \leq t_\ep$ will not affect the limit and the term inside the expectation is bounded, the rest of the proof can be completed in the same manner as the proof of Theorem \ref{th1B}.

We conclude this section with the

\mn
{\it Proof of Theorem \ref{th1C}.} Let $\mathcal{M}(t)$ be the three dimensional Poisson process defined in Section \ref{prelimsec}. Using (\ref{expident}), we see that in order for the contribution of $Z_1^{s,x}(t)$ to the limit of $t^{1+p} e^{-(\lambda_0+b)t} Z_1(t)$ to be $> z$ we need
$$
v > z t^{-(1+p)} e^{(b-x)t} e^{(\lambda_0+x)s}
$$
Therefore, the expected number of mutations that contribute more than $z$ to the limit is
$$
u_1 V_0 \int_0^b dx \, g(x) \int_0^t ds \, e^{\lambda_0 s} \, \frac{\lambda_0+x}{a_0+x}
\exp\left( - \frac{\lambda_0+x}{a_0+x} \cdot z t^{-(1+p)} e^{(b-x)t} e^{(\lambda_0+x)s} \right)
$$
In order to turn the big exponential into $e^{-r}$ we change variables:
$$
s = \frac{1}{\lambda_0+x} \log \left(\frac{r}{z t^{-(1+p)} e^{(b-x)t} \frac{\lambda_0+x}{a_0+x}}\right)
$$
$ds = dr/r(\lambda_0+x)$ to get
\begin{align*}
u_1 V_0 \int_0^b dx & \, g(x) z^{-\lambda_0/(\lambda_0+x)} \left(\frac{\lambda_0+x}{a_0+x}\right)^{x/(\lambda_0+x)}
\cdot t^{(1+p)\lambda_0/(\lambda_0+x)}\\
& e^{-(b-x)t\lambda_0/(\lambda_0+x)}  \int_{\alpha(x,t)}^{\beta(x,t)}
\frac{dr}{\lambda_0+x} r^{-x/(\lambda_0+x)} e^{-r}
\end{align*}
where $\alpha(x,t) = z t^{-(1+p)} e^{(b-x)t} (\lambda_0+x)/(a_0+x)$ and
$\beta(x,t) = \alpha(x,t) e^{(\lambda_0+x)t}$.
As in the previous proof, the main contribution comes from $x \in [b_t,b]$ so when we change variables
$y = (b-x)t$, $dx = -dy/t$, replace the $x$'s by $b$'s and use $1=(1+p)\lambda_0/(\lambda_0+b)$
we convert the above into
$$
g(b) z^{-\lambda_0/(\lambda_0+b)} \frac{u_1V_0}{\lambda_0+b} \left(\frac{\lambda_0+b}{a_0+b}\right)^{b/(\lambda_0+b)}
\int_0^\infty dy \, e^{-y\lambda_0/(\lambda_0+b)} \int_0^\infty  r^{-b/(\lambda_0+b)} e^{-r}\, dr
$$
Performing the integrals gives the result with
$$
A_1(\lambda_0, b)= g(b) \frac{1}{\lambda_0}  \left(\frac{\lambda_0+b}{a_0+b}\right)^{b/(\lambda_0+b)}
\Gamma(\lambda_0/(\lambda_0+b))
\eqno\square
$$

\section{Bounded distributions, ${\bf Z}_{\bf k}$} \label{bndZksec}

We now move on to the proofs of Theorems \ref{th2A} and \ref{th2B}. Recall that we have defined $p_k$ by the relation
$$
k+p_k = \sum_{j=0}^{k-1} \frac{\lambda_0+kb}{\lambda_0+jb}.
$$

\mn
{\it Proof of Theorem \ref{th2A}.} Breaking things down according to the times and the sizes of the mutational changes we have
\begin{align} \label{Zkmean}
EZ_k(t) = \int_0^b dx_1 \, g(x_1) & \cdots \int_0^b dx_k\, g(x_k)  \int_0^t ds_1 \cdots \int_{s_{k-1}}^t ds_k \nonumber \\
& V_0 e^{\lambda_0 s_1} u_1 e^{(\lambda_0+x_1) (s_2-s_1)} \cdots u_k e^{(\lambda_0+x_1+ \cdots + x_k) (t-s_k)} \nonumber \\
= \int_0^b dx_1 \, g(x_1) & \cdots \int_0^b dx_k\, g(x_k)  \int_0^t ds_1 \cdots \int_{s_{k-1}}^t ds_k \nonumber \\
& V_0 u_1 \dotsm u_k e^{\lambda_0t} e^{x_1(t-s_1)} \cdots e^{x_k(t-s_k)}.
\end{align}
The first step is to show

\begin{lemma} \label{kt1}
Let $b_t = b - (2k+1)(\log t)/t$. The contribution to $EZ_k(t)$ from points $(x_1,\ldots x_k)$ with some $x_i \le b_t$
is $o(t^{-2k}e^{(\lambda_0+kb)t})$.
\end{lemma}

\mn
{\it Proof.} (\ref{expineq}) implies that
$$
\int_{s_{j-1}}^t ds_j \, e^{(x_j+\cdots+x_k)(t-s_j)}
= \frac{ e^{(x_j+\cdots+x_k)(t-s_{j-1})} -1 }{x_j+\cdots+x_k} \le t e^{(x_j+\cdots+x_k)(t-s_{j-1})}.
$$
Applying this and working backwards in the above expression for $EZ_k(t)$, we get
$$
EZ_k(t) \le t^k V_0 u_1 \cdots u_k  \int_0^b dx_1 \, g(x_1) \cdots \int_0^b dx_k\, g(x_k) e^{(\lambda_0 + x_1+\cdots+x_k)t}
$$
and the desired result follows. \eopt

\ms
With the Lemma established, when we work backwards
$$
\int_{s_{j-1}}^t ds_j \, e^{(x_j+\cdots+x_k)(t-s_j)}
= \frac{ e^{(x_j+\cdots+x_k)(t-s_{j-1})} -1 }{x_j+\cdots+x_k} \sim \frac{e^{(x_j+\cdots+x_k)(t-s_{j-1})}}{(k-j+1)b}
$$
From this and induction, we see that the contribution from points $(x_1,\ldots x_k)$ with $x_i \in[b_t,b]$ for all
$i$ is
$$
\sim \frac{ V_0 u_1 \cdots u_k }{b^k k!}  g(b)^k \int_{b_t}^b dx_1 \, \cdots \int_{b_t}^b dx_k\, e^{(\lambda_0 + x_1+\cdots+x_k)t}
$$
Changing variables $y_i = t(b-x_i)$ the above
$$
\sim \frac{ V_0 u_1 \cdots u_k g(b)^k }{b^k t^k k!} e^{(\lambda_0 + kb)t}
$$
which proves the desired result. \eopt

In the proof of the last result, we showed that the dominant contribution comes from mutations with $x_i > b_t$. To prove our limit theorem we will also need a result regarding the times at which the mutations to the dominant types occur.

\begin{lemma} \label{kt2}
Let $\alpha_k = \frac{2k+1}{kb}$. The contribution to $EZ_k(t)$ from points with
$s_1 \ge \alpha_k \log t$ is $o(t^{-2k}e^{(\lambda_0+kb)t})$.
\end{lemma}

\begin{proof}
Replace the $x_i$'s in the exponents by $b$'s, we can see from \eqref{Zkmean} that the expected contribution from points with $s_1 \geq \alpha_k \log t$ is
\begin{align*}
& \leq b^k G^k V_0 u_1 \dotsm u_k \int_{\alpha_k \log t}^t ds_1 \int_{s_1}^t ds_2 \dotsm \int_{s_{k-1}}^t ds_k e^{\lambda_0 t} e^{b(t-s_1)} \dotsm e^{b(t-s_k)} \\
& \leq Ce^{\lambda_0t} \int_{\alpha_k \log t}^t e^{kb(t-s_1)} ds_1  \\
& \leq C e^{(\lambda_0+kb)t} t^{-\alpha_kkb}
\end{align*}
and the desired result follows.
\end{proof}

Recall that
$$
k+p_k = \sum_{j=0}^{k-1} \frac{\lambda_0+kb}{\lambda_0+jb}.
$$
For the induction used in the next proof, we will also need the corresponding quantity with $\lambda_0$ replaced by $\lambda_0+x$
and $k$ by $k-1$
$$
k-1+p_{k-1}(x) = \sum_{j=0}^{k-2} \frac{\lambda_0+x+(k-1)b}{\lambda_0+x+jb}
$$
which means
$$
p_{k-1}(x) = \sum_{j=0}^{k-2} \frac{(k-1-j)b}{\lambda_0+x+jb}
$$
The limit will depend on the mutation rates through
$$
u_{1,k} = \prod_{j=1}^k u_j^{\lambda_0/(\lambda_0+(j-1)b)}
$$
Again we will need the corresponding quantity with $k-1$ terms
$$
u_{2,k}(x) = \prod_{j=1}^{k-1} u_{j+1}^{(\lambda_0+x)/(\lambda_0+x+(j-1)b)}.
$$
We shall write $u_{2,k} = u_{2,k}(b)$ and note that
\beq
u_{1,k} = u_1 u_{2,k}^{\lambda_0/(\lambda_0+b)}
\label{urec}
\eeq

\mn
{\it Proof of Theorem \ref{th2B}.} We shall prove the result under the more general assumption that $Z_0(t) \sim V_0 e^{\lambda_0 t}$ for some constant $V_0$. The result then holds for $k=1$ by Lemma \ref{th1Bgen}. We shall prove the general result by induction on $k$. To this end, suppose the result holds for $k-1$. Let $Z_k^{s,x,v}(t)$ be the type-$k$ descendants at time $t$ of the 1 mutant at $(s,x,v) \in \mathcal{M}(t)$. Since $Z_1^{s,x}(t) \sim v e^{(\lambda_0+x)(t-s)}$ compared to $Z_0(t) \sim V_0 e^{\lambda_0t}$, it follows from the induction hypothesis that
\begin{align}
E\exp\bigl( - \theta (t-s)^{k-1+p_{k-1}(x)} & e^{-(\lambda_0+x+(k-1)b)(t-s)} Z_k^{s,x,v}(t) \bigr)
\nonumber\\
&\to \exp\bigl( - c_{k-1}(\lambda_0+x,b) v u_{2,k}(x) \theta^{(\lambda_0+x)/(\lambda_0+x+(k-1)b)} \bigr)
\label{philim2}
\end{align}
Integrating over the contributions from the three-dimensional point process we have
\begin{align*}
E\exp(-\theta Z_k(t) ) = \exp\biggl( -  &\int_0^b dx \, g(x) \int_0^t ds \, u_1 V_0 e^{\lambda_0 s} \\
& \int_0^\infty dv \, \left( \frac{\lambda_0+x}{a_0+x} \right)^2 \exp\left( - \frac{\lambda_0+x}{a_0+x} v \right)
( 1 - \phi^{k-1}_{x,v,t-s}(\theta)) \biggr)
\end{align*}
where $\phi^{k-1}_{x,v,t-s}(\theta) = E \exp(-\theta Z_k^{0,x,v}(t-s))$. To prove the desired result we need to
replace $\theta$ by $\theta t^{k+p_k} e^{-(\lambda_0+kb)t}$. Doing this with (\ref{philim2}) in mind
we have
\begin{align*}
E\exp(-\theta t^{k+p_k}  e^{-(\lambda_0+kb)t} & Z_k(t) )  \\
=  \exp\biggl( - \int_0^b  & dx \, g(x) \int_0^t  ds \, u_1 V_0 e^{\lambda_0 s}
\int_0^\infty dv \, \left( \frac{\lambda_0+x}{a_0+x} \right)^2
\exp\left( - \frac{\lambda_0+x}{a_0+x} v \right) \\
\nass
& \left\{ 1 - \phi^{k-1}_{x,v,t-s}( \theta t^{k+p_k}
e^{-(\lambda_0+x+(k-1)b)(t-s)} e^{-(b-x)t} e^{-(\lambda_0+x+b(k-1))s} ) \right\}
\biggr)
\end{align*}
By Lemmas \ref{kt1} and \ref{kt2}, we can restrict attention to $x \in [b_t,b]$ and $s \le \alpha_k \log t$. The first restriction
implies that all of the $x$'s except the one in $(b-x)$ can be set equal to $b$ and the second that we can replace $t$ by $t-s$.
Since $(k+p_k) - (k-1+p_{k-1}(b)) = (\lambda_0+kb)/\lambda_0$, the term in the exponential is
\begin{align*}
= - \int_{b_t}^b & dx \, g(x) \int_0^{\alpha_k\log t} ds \, u_1 V_0 e^{\lambda_0 s}
\int_0^\infty dv \, \left( \frac{\lambda_0+b}{a_0+b} \right)^2  \exp\left( - \frac{\lambda_0+b}{a_0+b} v \right) \\
\nass
& \left( 1 - \phi_{x,v,t-s}( \theta (t-s)^{k-1+p_{k-1}(b)} e^{-(\lambda_0+kb)(t-s)}
t^{(\lambda_0+kb)/\lambda_0} e^{-(b-x)t} e^{-(\lambda_0+kb)s} ) \right)
\end{align*}
Changing variables $s=R(t)+r$ where $R(t)=(1/\lambda_0)(\log t)$, and $y=(b-x)t$, $dy = -tdx$  the above becomes
\begin{align*}
& = - g(b) \int_0^{(2k+1)\log t} dy
\int_0^\infty dv \, \left( \frac{\lambda_0+b}{a_0+b} \right)^2  \exp\left( - \frac{\lambda_0+b}{a_0+b} v \right) \\
\nass
& \qquad \int_{-R(t)}^{\alpha_k\log t-R(t)} dr \, u_1 V_0 e^{\lambda_0r}
\left( 1 - \phi^{k-1}_{x,v,t-s}( \theta (t-s)^{k-1+p_{k-1}(b)} e^{-(\lambda_0+kb)(t-s)} e^{-y} e^{-(\lambda_0+kb)r} ) \right)
\end{align*}
Using (\ref{philim2}) now we have that the $1-\phi$ term converges to
$$
1 - \exp\left(-c_{k-1}(\lambda_0+b,b) v u_{2,k} [\theta e^{-y}] ^{(\lambda_0+b)/(\lambda_0+kb)} e^{-(\lambda_0+b)r}\right)
$$
To simplify the exponential we let
$$
r = \frac{1}{\lambda_0+b} ( q + Q(v,y))
\quad\hbox{where}\quad Q(v,y)= \log\left\{ c_{k-1}(\lambda_0+b,b) v u_{2,k} [\theta e^{-y}]^{(\lambda_0+b)/(\lambda_0+kb)}
\right\}
$$
$dr = dq/(\lambda_0+b)$. Plugging this into $e^{\lambda_0r}$ results in
$$
e^{q\lambda_0/(\lambda_0+b)} (c_{k-1}(\lambda_0+b,b) v u_{2,k})^{\lambda_0/(\lambda_0+b)} \theta^{\lambda_0/(\lambda_0+kb)}
e^{-y\lambda_0/(\lambda_0+kb)}
$$
so the exponential converges to
\begin{align*}
-g(b) \frac{c_{k-1}(\lambda_0+b,b)^{\lambda_0/(\lambda_0+b)}}{\lambda_0+b}
 V_0 u_1 u_{2,k}^{\lambda_0/(\lambda_0+b)}\theta^{\lambda_0/(\lambda_0+kb)} & \\
\int_0^\infty dv \, \left( \frac{\lambda_0+b}{a_0+b} \right)^2
v^{\lambda_0/(\lambda_0+b)} & \exp\left( - \frac{\lambda_0+b}{a_0+b} v \right) \\
 \int_0^\infty dy \, e^{-y\lambda_0/(\lambda_0+kb)}
\int_{-\infty}^{\infty} & \frac{dq}{\lambda_0+b} \, e^{q\lambda_0/(\lambda_0+b)} (1 - \exp(- e^{-q} ))
\end{align*}
To clean this up, we note that letting $w = v(\lambda_0+b)/(a_0+b)$, $dw = dv(\lambda_0+b)/(a_0+b)$
\begin{align}
\int_0^\infty dv \, \left( \frac{\lambda_0+b}{a_0+b} \right)^2  &
v^{\lambda_0/(\lambda_0+b)} \exp\left( - \frac{\lambda_0+b}{a_0+b} v \right)
\nonumber\\
& = \left( \frac{a_0+b}{\lambda_0+b}\right)^{-1 + \lambda_0/(\lambda_0+b)}
\Gamma(1+\lambda_0/(\lambda_0+b))
\label{ci1}
\end{align}
The second integral is easy:
\beq
\int_0^\infty dy \, e^{-y\lambda_0/(\lambda_0+kb)} = \frac{\lambda_0+kb}{\lambda_0}
\label{ci2}
\eeq
The third one looks weird but when you put $x= e^{-q}$, $dx = -e^{-q}\, dq$, or $dq = - dx/x$ it is
$$
= \int_0^\infty dx \, x^{-1-\lambda_0/(\lambda_0+b)} (1-e^{-x}) \, dx
$$
then integrating by parts $f(x)=1-e^{-x}$, $g'(x) = x^{-1-\lambda_0/(\lambda_0+b)}$, $f'(x) = e^{-x}$,
$g(x) = x^{-\lambda_0/(\lambda_0+b)}(\lambda_0+b)/\lambda_0$
turns it into
\beq
\frac{\lambda_0+b}{\lambda_0} \Gamma(1 - \lambda_0/(\lambda_0+b))
\label{ci3}
\eeq
Putting this all together and using \eqref{urec}, we have
\begin{align*}
& c_{k-1}(\lambda_0+b,b)^{\lambda_0/(\lambda_0+b)} \cdot g(b) \frac{\lambda_0+kb}{\lambda_0} \cdot
 V_0 u_{1,k} \theta^{\lambda_0/(\lambda_0+kb)} \\
& \qquad \cdot \frac{1}{\lambda_0} \left( \frac{a_0+b}{\lambda_0+b}\right)^{-1 + \lambda_0/(\lambda_0+b)}
\Gamma(1+\lambda_0/(\lambda_0+b)) \Gamma(1 - \lambda_0/(\lambda_0+b))
\end{align*}
Setting $c_k(\lambda_0,b)$ equal to the quantity in the last display
divided by $V_0 u_{1,k} \theta^{\lambda_0/(\lambda_0+kb)}$ we have proved the result. \eopt

\medskip
To work out an explicit formula for the constant and to compare with Durrett and Moseley (2009), it is useful to
let $\lambda_j = \lambda_0+jb$, $a_j = a_0+jb$ and
$$
c_{h,j} = \frac{1}{\lambda_{j-1}} \left( \frac{a_j}{\lambda_j}\right)^{-1 + \lambda_{j-1}/\lambda_j}
\Gamma(1+\lambda_{j-1}/\lambda_j) \Gamma(1 - \lambda_{j-1}/\lambda_j)
$$
From this we see that
\begin{align*}
c_k(\lambda_0,b)  =  &c_{k-1}(\lambda_1,b)^{\lambda_0/\lambda_1} g(b) \frac{\lambda_k}{\lambda_0} c_{h,1} \\
=  & c_{k-2}(\lambda_2,b)^{\lambda_0/\lambda_2}
\cdot \left( g(b) \frac{\lambda_{k-1}}{\lambda_0} c_{h,2} \right)^{\lambda_0/\lambda_1} \cdot g(b) \frac{\lambda_k}{\lambda_0} c_{h,1}
\end{align*}
and hence
$$
c_k(\lambda_0,b) = \prod_{j=1}^{k} \left( g(b) \frac{\lambda_{k-j+1}}{\lambda_0} c_{h,j} \right)^{\lambda_0/\lambda_{j-1}}
$$

In Durrett and Moseley (2009) if we let ${\cal F}_{k-1}$ be the $\sigma$-field generated by $Z_j(t)$ for $j\le k$ and all $t\ge 0$ then
$$
E(e^{-\theta V_k}|{\cal F}_{k-1}) = \exp( - u_k V_{k-1}c_{h,k} \theta^{\lambda_{k-1}/\lambda_k})
$$
Iterating we have
\begin{align*}
E(e^{-\theta V_k}|{\cal F}_{k-2}) & = E( \exp( - u_k V_{k-1}c_{h,k}\theta^{\lambda_{k-1}/\lambda_k}) |{\cal F}_{k-2}) \\
& = \exp\left( -u_{k-1} u_k^{\lambda_{k-2}/\lambda_{k-1}} V_{k-2}
c_{h,k-1} c_{h,k}^{\lambda_{k-2}/\lambda_{k-1}}  \theta^{\lambda_{k-2}/\lambda_k}  \right)
\end{align*}
and hence
$$
E(e^{-\theta V_k}|V_0) = \exp( - c_{\theta,k} V_0 u_{1,k} \theta^{\lambda_0/\lambda_k} )
$$
where $c_{\theta,k} = \prod_{j=1}^{k} c_{h,j}^{\lambda_0/\lambda_{j-1}}$.

\section{Proofs for unbounded distributions} \label{unbndsec}

In this Section, we prove Theorem \ref{uthW}. The first step is to show that unlike in the case of bounded mutational advances, for unbounded distributions, the main contribution to the limit is given by the descendants of a single mutations. The largest growth rate will come from $z = O(t^{(\alpha+1)/\alpha})$ so the next result is enough. Recall that the mean number of mutations with growth rate larger than $z$ has
\begin{align*}
\mu(z,\infty) & = V_0 u_1 \int_0^t \left( \frac{z}{t-s} - \lambda_0 \right)^\beta e^{\lambda_0s}
\exp\left( - \gamma \left( \frac{z}{t-s} - \lambda_0 \right)^{\sqz\alpha} \, \right) \, ds \\
&= V_0 u_1 \int_0^t \left( \frac{z}{t-s} - \lambda_0 \right)^\beta \exp(\phi(s,z)) \, ds
\end{align*}
where $\phi$ is as in \eqref{maxme}.

\begin{lemma} \label{smallneg}
Let $\bar{z} > \lambda_0 t$. Then
\beqx
E \left(\sum_{(s,x):z(s,x) \leq \bar{z}} Z_1^{s,x}(t) \right) \leq C V_0 u_1 \bar{z} e^{\lambda_0 t + \bar{z}}
\eeqx
as $t \to \infty$.
\end{lemma}

\begin{proof}
The expected number of individuals produced by mutations with growth rates $ \le \bar{z}$ is
$$
V_0 u_1 \int_0^t \int_0^{\frac{\bar{z}}{t-s} - \lambda_0} e^{\lambda_0 s} \cdot y^\beta e^{-\gamma y^\alpha} \cdot e^{z(s,y)} \, dy \, ds.
$$
Changing variables $y \mapsto u = z(s,y)$, that is $y= u/(t-s) - \lambda_0$, $dy = du/(t-s)$, and using Fubini's theorem to switch the order of integration, we can see that the above is
\beq
\le V_0 u_1 e^{\lambda_0 t + \bar{z}} \int_{0}^{\bar{z}}  \int_0^t (u/(t-s)-\lambda_0)^{\beta} \exp\left( - \gamma \left( \frac{u}{t-s} -\lambda_0\right)^\alpha \right) \, \frac{ds}{(t-s)} \, du.
\label{int2}
\eeq
But then if we change variables $s \mapsto r= u/(t-s)-\lambda_0$, $dr = u ds /(t-s)^2$, we can see that the inner integral is
\beqx
\le \int_{-\lambda_0}^\infty \frac{r^\beta}{r+\lambda_0}e^{-\gamma r^\alpha} dr \le C
\eeqx
yielding the desired bound.
\end{proof}

To motivate the proof of the general result, we begin with the case when $\alpha = 1$.

\begin{proof} [Proof of Theorme \ref{uthexp}.] Since
\beqx
Z_1(t) = \sum_{(s,x) \in \calN(t)} Z_1^{s,x}(t) = \sum_{(s,x): z(s,x) \leq z} Z_1^{s,x}(t) + \sum_{(s,x): z(s,x) > z} Z_1^{s,x}(t)
\eeqx
for any $z >0$, we have
\beqx
\frac{1}{t} \log Z_1(t) \sim \frac{1}{t} \left[ \log \left(\sum_{(s,x): z(s,x) \leq z} Z_1^{s,x}(t)\right) \vee  \log \left(\sum_{(s,x): z(s,x) > z} Z_1^{s,x}(t)\right) \right]
\eeqx
as $t \to \infty$. Lemma \ref{smallneg} tells us that if there is a mutation with growth rate $z = O(t^2)$, then the contribution from mutations with growth rates smaller than $z-\ep$ can be ignored so it suffices to describe the distribution of the largest growth rates. We will show that
\beq \label{goal}
\mu(z,\infty) \to \begin{cases} 4 c_0^{\beta} ( \pi/\lambda_0)^{1/2} V_0 u_1 \exp(\gamma \lambda_0- 2 \lambda_0x/2c_0) & \text{ if } \quad z=c_0t^2\left(1 + \frac{(2\beta+1) \log t}{\lambda_0 t} + \frac{x}{c_0t } \right) \\
0 & \text{ if } \quad z \gg c_0t^2\left(1 + \frac{(2\beta+1) \log t}{\lambda_0 t}\right)\end{cases}
\eeq
so that the largest growth rate is $O(t^2)$ and comes from the rightmost particle in the point process with intensity given by \eqref{ptp}.

To prove \eqref{goal}, we first need to locate the maximum of $\phi$. Let $z > \lambda_0 t$ so that there exists a unique maximum $s_z$. Solving $\phi_s(s,z) = 0$ and using the expression for $\phi_s$ in \eqref{der1} yields
\beqx
s_z = t - a_0 z^{1/2}
\eeqx
where $a_0 = (\gamma/\lambda_0)^{1/2} = (4c_0)^{-1/2}$ which leads to the expression
\begin{align} \label{phimax}
\phi(s_z,z) &= \lambda_0 t - \lambda_0(t-s_z) - \gamma \left(\frac{z}{t-s} -\lambda_0\right) \nonumber \\
& = \lambda_0 t - \lambda_0 a_0 z^{1/2} - \gamma z^{1/2}/a_0 + \gamma \lambda_0 \nonumber \\
&= \lambda_0 (t - 2a_0 z^{1/2}) + \gamma \lambda_0.
\end{align}
If we take
\beqx
z_x =c_0t^2\left(1 + \frac{\kappa \log t}{t} + \frac{x}{c_0t } \right) = \left(\frac{t}{2a_0}\right)^2 \left(1 + \frac{\kappa \log t}{t} + \frac{4 a_0^2 x}{t} \right)
\eeqx
in \eqref{phimax} and use $(1+y)^{1/2} = 1 + y/2 + O(y^2)$, we obtain
\begin{align} \label{phimaxmax}
\phi(s_{z_x},z_x) &=  - \frac{\lambda_0 \kappa \log t}{2} - 2 \lambda_0 a_0^2 x  + \gamma \lambda_0 + o(1)
\end{align}
as $t \to \infty$. Furthermore, \eqref{der2} implies that
\begin{align*}
\phi_{ss}(s_{z_x},z_x) = -\frac{2\gamma z_x}{(t-s_{z_x})^3} =  -\frac{2 \gamma}{a_0^3 z_x^{1/2}} &= -\frac{a}{t} + o(1) \\
\phi_{sss}(s_{z_x},z_x) = -\frac{6\gamma z_x}{(t-s_{z_x})^4}  = -\frac{6 \gamma}{a_0^4 z} &= -\frac{24 \gamma}{a_0^2 t^2} + o(1)
\end{align*}
as $t \to \infty$ with $a= 4 \gamma/a_0^2$. Since $\phi_s(s_z,z) = 0$, taking a Taylor expansion around $s_z$ yields
\beq \label{phitaylor}
\phi(s,z_x) = \phi(s_{z_x},z_x) -\frac{a}{2 t} (s-s_{z_x})^2 + g(s,z_x)
\eeq
where $|g(s,z)| \leq C |s-s_{z}|^3/t^2$ for all $s$. Also note that letting
\beqx
\psi(s,z) = \left(\frac{z}{t-s} - \lambda_0\right)^\beta
\eeqx
we have
\begin{align*}
\psi(s_{z_x},z_x) &= \left(\frac{z_x}{t-s_{z_x}} - \lambda_0\right)^\beta \\
&= z_x^{\beta/2}/a_0^\beta + o(z_x^{\beta/2}) \\
&= (2c_0)^\beta  t^\beta + o(t^\beta)
\end{align*}
so that
\begin{align*}
\psi(s,z_x) = (2c_0)^\beta  t^\beta + g_2(s,z_x)
\end{align*}
where $|g_2(s,z)||s-s_z|^{-1}t^{-\beta} = o(1)$.

Write
\begin{align*}
 \int_0^t \psi(s,z_x) e^{\phi(s,z_x)} \, ds = \int_A  \psi(s,z_x) e^{\phi(s,z_x)} \, ds + \int_{A^c}\psi(s,z_x) e^{\phi(s,z_x)} \, ds
\end{align*}
where $A = \{s: |s-s_{z_x}| \leq C (t\log t)^{1/2}\} \cap [0,t]$. Since concavity implies that for $s \in A^c$ and $C$ sufficiently large, we have
$$
\exp(\phi(s,z_x)) \le \frac{1}{t^{2+\beta}} \exp(\phi(s_{z_x},z_x))
$$
the contribution of the second integral is negligible. After the change of variables $s= s_{z_x} + (t/a)^{1/2}r$, when $t$ is large, the first integral becomes
\begin{align*}
 \int_A \psi(s,z_x)  e^{\phi(s,z_x)} \, ds = ((2c_0)^\beta  t^\beta + o(1)) e^{\phi(s_{z_x},z_x)} \int_{-C (\log t)^{1/2}}^{C (\log t)^{1/2}} e^{g(s,z_x)} e^{-r^2/2} (t/a)^{1/2} \, dr.
\end{align*}
and therefore since $|g(s,z_x)| \leq  C (t \log t)^{3/2}/t^2 $ when $s \in A$, we have
\beq \label{intensity}
\mu(z_x,\infty) = V_0 u_1 \int_0^t \psi(s,z_x)  e^{\phi(s,z_x)} \, ds \sim  b V_0 u_1 t^{\beta+ 1/2} e^{\phi(s_{z_x},z_x)}
\eeq
where $b = (2c_0)^\beta \sqrt{2\pi/a} = (2c_0)^\beta  ( \pi/\lambda_0)^{1/2}$. Since
\begin{align*}
\phi(s_{z_x},z_x) &=  - \frac{\kappa \lambda_0 \log t}{2} - 2 \lambda_0 a_0^2 x  + \gamma \lambda_0
\end{align*}
we can conclude that
$$
\mu(z_x,\infty) \to \begin{cases} V_0 u_1 b \exp(\gamma \lambda_0- 2 \lambda_0 a_0^2 x) = V_0 u_1 b \exp(\gamma \lambda_0- 2 \lambda_0x/2c_0)  & \text { if } \quad \kappa = \frac{2\beta +1}{\lambda_0} \\ 0 & \text{ if } \quad \kappa > \frac{2\beta +1}{\lambda_0} \end{cases}
$$
which proves \eqref{goal} since this argument remains true even if $\kappa= \kappa(t)$ and $\liminf \kappa(t) > \frac{2\beta+1}{\lambda_0}$.
\end{proof}

When $\alpha \neq 1$, we no longer have an explicit formula for the maximum value $s_z$ which complicates the process of identifying the largest growth rate. We shall assume for convenience that $\alpha >0$ is an integer.

\mn
{\it Proof of Theorem \ref{uthW}.} As in the proof of Theorem \ref{uthexp}, it suffices to describe the distribution for the largest growth rates. Let $z > \lambda_0 t$ so the maximum $s_z$ exists. To find a useful expression for the value of $\phi(s_z,z)$, we write
$$
\phi(s,z) = \lambda_0 t - \lambda_0(t-s) - \gamma \left( \frac{z}{t-s} - \lambda_0 \right)^{\sqz\alpha}.
$$
Using the definition of $s_z$ as the solution to $\phi_s(s_z,z)= 0$ yields the condition that
\begin{align*}
(t-s_z)^{\alpha+1} = \frac{\alpha \gamma}{\lambda_0} z^\alpha (1 - \lambda_0 \frac{t-s_z}{z})^{\alpha-1}
\end{align*}
i.e.,
\begin{align*}
t-s_z = \left(\frac{\alpha \gamma}{\lambda_0}\right)^{1/(\alpha+1)} z^{\alpha/(\alpha+1)} \left(1 - \lambda_0 \frac{t-s_z}{z}\right)^{(\alpha-1)/(\alpha+1)}.
\end{align*}
If we substitute the right side of this equation back in for $t-s_z$ in the parenthesis, then writing $a_0 = (\alpha\gamma/\lambda_0)^{1/(\alpha+1)}$, we have
\begin{align*}
t-s_z &= a_0 z^{\alpha/(\alpha+1)} \left(1 - \lambda_0 a_0 z^{-1/(\alpha+1)}\left(1-\frac{\lambda_0(t-s_z)}{z}\right)^{\frac{\alpha-1}{\alpha+1}}\right)^{\frac{\alpha-1}{\alpha+1}}  \\
& = a_0 z^{\alpha/(\alpha+1)} \left(1 - \lambda_0 a_0 z^{-1/(\alpha+1)}\left(1 - \lambda_0 a_0 z^{-1/(\alpha+1)}\left(1- \frac{\lambda_0(t-s_z)}{z}\right)^{\frac{\alpha-1}{\alpha+1}}\right)^{\frac{\alpha-1}{\alpha+1}}  \right)^{\frac{\alpha-1}{\alpha+1}}
\end{align*}
\noindent
We repeat this $\alpha$ times and then use the approximation $(1-x)^n = 1 - nx+ O(x^2)$ repeatedly with  $n = (\alpha-1)/(\alpha+1)$ to obtain
\begin{align} \label{diff}
t-s_z = z^{\alpha/(\alpha+1)} \left(\sum_{j=0}^\alpha a_j z^{-j/(\alpha+1)} + O(z^{-1}) \right)
\end{align}
where
$$
a_j = a_0 \left(\frac{\lambda_0 a_0 (\alpha-1)}{\alpha+1}\right)^j
$$
for $j \geq 1$. The error term is $O(z^{-1})$ because
$$
 0 < (1 - \lambda_0 (t-s)/z) \leq 1
$$
for all $z > \lambda_0 t $ and $s \leq t$. Factoring out $a_0$ in \eqref{diff} and using $(1+x)^{-1} = \sum (-x)^j$ when $|x| < 1$, we have that
\begin{align} \label{comterm}
\frac{z}{t-s}  - \lambda_0 &=  a_0^{-1} z^{1/\alpha+1}\left(1 - \sum_{i_1=1}^{\alpha} a_0^{-1} a_{i_1} z^{-i_1/(\alpha+1)} + \sum_{i_1,i_2=1}^\alpha  a_0^{-2} a_{i_1}a_{i_2} z^{-(i_1+i_2)/(\alpha+1)} \right. \nonumber \\
&\left.  - \dotsm + (-1)^\alpha \sum_{i_1,...,i_\alpha = 1}^\alpha a_0^{-\alpha} \prod_{j=1}^\alpha a_{i_j} z^{-\sum_{j=1}^\alpha i_j/(\alpha+1)}  + O(z^{-1}) \right) - \lambda_0 z^{1/(\alpha+1)}z^{-1/(\alpha+1)} \nonumber \\
& =  z^{1/(\alpha+1)} \left( \sum_{j=0}^\alpha b_j z^{-j/(\alpha+1)} + O(z^{-1}) \right)
\end{align}
for large $z$ where the $b_j$ are given by
\begin{align*}
b_0 &= 1/a_0 \\
b_1 &= -a_1/a_0^2 - \lambda_0 \\
b_2 &= -(a_2-a_1^2)/a_0^3 \\
b_3 &= -(a_4 - 2a_1a_3 - a_2^2 -3a_1^2 a_2 + a_1^4)/a_0^4
\end{align*}
and in general,
\beqx
b_i = \sum_{k=1}^\alpha \sum_{i_1,...,i_k: i_1 + \dotsm +i_k = i} (-a_0)^{-(k+1)} \prod_{j=1}^k a_{i_j}.
\eeqx
\eqref{comterm} implies that
\begin{align*}
-\gamma \left(\frac{z}{t-s} - \lambda_0\right)^\alpha &= -\gamma z^{\alpha/(\alpha+1)} \left( b_0^\alpha + \alpha b_0^{\alpha-1}b_1 z^{-1/(\alpha+1)} \right. \\
& \quad \left. + \left(\alpha b_0^{\alpha-1}b_2 + \binom{\alpha}{2}b_0^{\alpha-2}b_1^2 \right)z^{-2/(\alpha+1)} + \dotsm \right. \\
& \quad + \left. \left(\alpha b_0^{\alpha-1}b_\alpha + \dotsm + b_1^\alpha\right) z^{\alpha/(\alpha+1)} + O(z^{-1})\right)
\end{align*}
and therefore,
\begin{align} \label{phimax2}
\phi(s_z,z) &= \lambda_0 t + \lambda_0(t-s) - \gamma \left( \frac{z}{t-s} - \lambda_0 \right)^{\sqz\alpha}  \nonumber\\
&= \lambda_0t + \sum_{j =0}^{\alpha} d_j z^{\frac{\alpha-j}{\alpha+1}} + O(z^{-1/(\alpha+1)})
\end{align}
where the $d_j$ can be calculated explicitly, for example:
\begin{align*}
d_0 &= -\lambda_0 a_0 -\gamma b_0^\alpha \\
d_1 &= -\lambda_0 a_1 - \gamma\alpha b_0^{\alpha-1}c_1 \\
d_2 &= -\lambda_0 a_2 - \gamma \left(\alpha b_0^{\alpha-1}b_2 + \binom{\alpha}{2}b_0^{\alpha-2}b_1^2 \right) \\
d_3 &= -\lambda a_3 - \gamma \left(\alpha b_0^{\alpha-1}b_3 + \binom{\alpha}{2} b_0^{\alpha-2}b_1 b_2 + \binom{\alpha}{3}b_0^{\alpha-3}b_1^3  \right).
\end{align*}
To figure out the distribution of the growth rate for the largest mutant, we let $c_0 = (-\lambda_0/d_0)^{(\alpha+1)/\alpha}$ and then search for $\kappa_j$, $j =1,...,\alpha-1$ and $\kappa$ so that plugging
\beqx
z_x = c_0 t^{(\alpha+1)/\alpha}\left(1+ \sum_{j=1}^{\alpha-1} \kappa_j t^{-j/\alpha} + \frac{x}{c_0t} + \frac{\kappa \log t}{t}  \right)
\eeqx
into \eqref{phimax2} yields
\beq \label{phimaxmax2}
\phi(s_{z_x},z_x) = k_1 - k_2 x - k_3\log t
\eeq
for some constants $k_1$, $k_2$, $k_3$. Substituting $z_x$ into \eqref{phimax2} and writing $\kappa_0 = 1$, $\kappa_\alpha = x/c_0 $ to ease the notation we obtain
\begin{align*}
\phi(s_{z_x},z_x) &= \lambda_0t + \sum_{j=0}^{\alpha} d_j \left(-\frac{\lambda_0t}{d_0}\right)^{(\alpha-j)/\alpha}\left(\sum_{j=0}^{\alpha} \kappa_j t^{-j/\alpha} + \kappa t^{-1} \log t \right)^{(\alpha-j)/(\alpha+1)} + O(t^{-1/\alpha}).
\end{align*}
Since $\lambda_0 t + d_0(-\lambda_0 t/d_0) = 0$, the first order terms in this expansion is $t^{(\alpha-1)/\alpha}$ and after using the Taylor series expansion
$$
(1+x)^p = 1 + p x + p(p-1)x^2/2  + \dotsm + p(p-1)\dotsm (p-\alpha+1)x^\alpha/\alpha! + O(x^{\alpha+1})
$$
we obtain
\begin{align} \label{phiz0}
\phi(s_{z_0},z_0) &= \sum_{j=1}^{\alpha} \rho_j t^{(\alpha-j)/\alpha}  + \rho \log t+ O(t^{-1/\alpha}\log t)
\end{align}
where
\begin{align*}
\rho &=  d_0 \left(-\frac{\lambda_0}{d_0}\right) \left(\frac{\alpha}{\alpha+1}\right)\kappa = - \frac{\alpha\lambda_0}{\alpha+1} \kappa \\
\rho_1 &= d_0 \left(-\frac{\lambda_0}{d_0}\right) \left(\frac{\alpha}{\alpha+1}\right)  c_1 + d_1\left(-\frac{\lambda_0}{d_0}\right)^{(\alpha-1)/\alpha} \\
\rho_2 &=  d_0 \left(-\frac{\lambda_0}{d_0}\right) \left[ \frac{\alpha}{\alpha+1} c_2 +  \frac{\alpha}{\alpha+1}\left(\frac{\alpha}{\alpha+1}-1\right) c_1^2\right] \\
& \quad + d_1 \left(-\frac{\lambda_0}{d_0}\right)^{(\alpha-1)/\alpha} \left(\frac{\alpha-1}{\alpha}\right) c_1 + d_2 \left(-\frac{\lambda_0}{d_0}\right)^{(\alpha-2)/\alpha}
\end{align*}
and in general
\begin{align*}
\rho_j &= \sum_{i=0}^j d_i\left(-\frac{\lambda_0}{d_0}\right)^{(\alpha-i)/\alpha} \sum_{k=1}^{j-i} \prod_{\ell=1}^{k}\left(\frac{\alpha-i}{\alpha+1}-\ell+1\right) \kappa_{i_\ell}
\end{align*}
$j=1,1,...,\alpha$ where for each $i$ and $k$, in the inner product, $i_1,...,i_k$ are always chosen to satisfy $i_1 + i_2 + \dotsm + i_k = j-i$. Since $\rho_j$ depends only on $\kappa_i$, $i \leq j$, then after noting that the coefficient of $\kappa_j$ in $\rho_j$ is $-\alpha\lambda_0/(\alpha+1)$, we can use forward substitution to solve the system $\rho_j= 0$, $j=1,2,...,\alpha-1$ for $\kappa_j$ to obtain the recursive formulas
\beq \label{cis}
c_j \equiv \kappa_j =  - \frac{\alpha+1}{\alpha\lambda_0} \left(\rho_j  -   \frac{-\alpha\lambda_0}{\alpha+1} \kappa_j\right)
\eeq
for $i=1,2,...,\alpha-1$. Setting $\rho=-k_3$ yields
$$
\kappa= \frac{(\alpha+1)k_3}{\alpha \lambda_0}
$$
and for this choice of $c_j$, $\kappa$, we obtain \eqref{phimaxmax2} with
$$
k_2 = - \frac{\alpha}{\alpha+1}\frac{d_0}{c_0} \left(- \frac{\lambda_0}{d_0}\right) = \frac{\alpha \lambda_0}{(\alpha+1)c_0}
$$
and $k_1 = - (\rho_\alpha - k_2x).$ Since
\begin{align*}
\left(\frac{z_x}{t-s_{z_x}} - \lambda_0\right)^\beta &= z_x^{\beta/(\alpha+1)} /a_0^\beta + o(z_x^{\beta/(\alpha+1)}) \\
&= \left(\frac{c_0^{1/(\alpha+1)}}{a_0}\right)^\beta t^{\beta/\alpha} + o(z_x^{\beta/(\alpha+1)})
\end{align*}
choosing $k_3 = (2\beta/\alpha+1)/2$ replaces \eqref{phimaxmax} in the proof of Theorem \ref{uthexp}.

Now substituting \eqref{diff} and \eqref{comterm} in \eqref{der2} yields
\begin{align*}
\phi_{ss}(s_z,z) &= -  \alpha(\alpha-1) \gamma z^{\frac{\alpha-2}{\alpha+1}}\left(\sum_{j=0}^\alpha b_j z^{-j/(\alpha+1)}+ O(z^{-1})\right)^{\sqz\alpha-2}  \\
& \quad \times  \frac{z^2}{z^{4\alpha/(\alpha+1)} \left(\sum_{j=0}^\alpha a_j z^{-j/(\alpha+1)} + O(z^{-1})\right)^4}\\
& \quad -  \alpha \gamma z^{\frac{\alpha-1}{\alpha+1}}  \left(\sum_{j=0}^\alpha b_j z^{-j/(\alpha+1)} + O(z^{-1})\right)^{\sqz\alpha-1} \\
& \quad \times \frac{2z}{z^{3\alpha/(\alpha+1)}\left(\sum_{j=0}^\alpha a_j z^{-j/(\alpha+1)} + O(z^{-1})\right)^3} \\
&= [-\alpha(\alpha-1) \gamma b_0^{\alpha-2}/a_0^4 - \alpha\gamma b_0^{\alpha-1}/a_0^3] z^{-\alpha/(\alpha+1)} + o(z^{-\alpha/(\alpha+1)}) \\
&= - \frac{\alpha^2 \gamma}{a_0^{\alpha+2}}z^{-\alpha/(\alpha+1)} + o(z^{-\alpha/(\alpha+1)})
\end{align*}
where in the second to last line we have used the fact that $b_0 = a_0^{-1}$. When $z = z_x$, this becomes
\beqx
\phi_{ss}(s_{z_x},z_x) = - \frac{a}{t} + o(t^{-1})
\eeqx
where
$$
a= \frac{\alpha^2 \gamma}{a_0^{\alpha+2}c_0^{\alpha/(\alpha+1)}}.
$$
Since $\phi_s(s_z,z) = 0$ and a calculation similar to the one above shows that $\phi_{sss}(s_{z_x},z_x) = O(t^{-2})$, we have
$$
\phi(s,z_x) = \phi(s_{z_x},z_x) - \frac{a}{2}(s-s_{z_x})^2 + g(s_{z_x},z_x)
$$
where $|g(s,z)| \leq C |s-s_{z}|^3/t^2$ for all $s$. This replaces \eqref{phitaylor} from the $\alpha =1$ proof and the rest of the proof is the same. Note that the intensity for the limiting point process is given by
\beq \label{intensity2}
\left(\frac{c_0^{1/(\alpha+1)}}{a_0}\right)^\beta \sqrt{2\pi/a} \exp(k_1-k_2 x).
\eeq

 \eopt

\begin{remark} \label{timeofb}
From \eqref{diff}, we have
$$
t - s_{z_x} \sim  a_0 (c_0 t^{(\alpha+1)/\alpha})^{(\alpha+1)/\alpha} = \frac{\alpha t}{\alpha+1}
$$
which tells us that the time at which the mutant with largest growth rate is born is $\sim t/(\alpha+1)$.

\end{remark}

\section{Discussion} \label{discsec}
In this paper, we have analyzed a multi-type branching process model of tumor progression in which mutations increase the birth rates of cells by a random amount. We studied both bounded and unbounded distributions for the random fitness advances and calculated the asymptotic rate of expansion for the $k$th generation of mutants.

In the bounded setting, we found that there are only two parameters of the distribution that affect the limiting growth rate of the $k$th generation (see Theorems \ref{th1A}, \ref{th1B}, \ref{th2A}, and \ref{th2B}): the upper bound for the support of the distribution and the value of its density at the upper bound. This is a rather intuitive result since one would expect that in the long run, the $k$th generation will be dominated by mutants with the maximum possible fitness. In addition, we found that there is a polynomial correction to the exponential growth of the $k$th generation. This correction is not present in the case where the fitness advances are deterministic. We have discussed this point in further detail in Section 1.1 and after the proof of Theorem \ref{th2B} in Section \ref{bndZksec}. Finally, we showed that the limiting population is descended from several different mutations (see Theorem \ref{th1C}).

In the unbounded setting, we assumed that the distribution of the fitness advance has the form
$$
P(X>x)=x^{\beta}e^{-\gamma x^{\alpha}}
$$
where $\alpha, \beta$, and $\gamma$ are parameters. We found that the population of cells with a single mutation grows asymptotically at a super-exponential rate $\exp(t^{(\alpha+1)/\alpha})$ (see Theorems \ref{uthexp} and \ref{uthW}) and at large times, most of the first generation is derived from a single mutation (see Lemma \ref{smallneg}). The super-exponential growth rate suggests that the exponential distribution, which is often used for the fitness advances of an organism due to natural selection, is not a good choice for modeling the mutational advances in the progression to cancer where there is very little evidence for populations growing at a super-exponential rate.


These conclusions provide several interesting contributions to the existing literature on evolutionary models of cancer progression. First, our model generalizes previous multi-type branching models of tumor progression by allowing for random fitness advances as mutations are accumulated and provides a mathematical framework for further investigations into the role played by the fitness distribution of mutational advances in driving tumorigenesis. Second, we have discovered that bounded distributions lead to exponential growth whereas unbounded distributions lead to super-exponential growth. This dichotomy might provide a new method for testing whether a tumor population has evolved with an unbounded distribution of mutational advances. Third, we observe that in the case of bounded distributions, the growth rate of the tumor is somewhat `robust' with respect to the mutational fitness distribution and depends only on its upper endpoint. Finally, our calculations of the growth rates for the $k$th generation of mutants serve as a groundwork for studying the evolution and role of heterogeneity in tumorigenesis. These implications will be explored further in future work.


\section*{References}

\frenchspacing

\mn
Becskei A., Kaufmann B.B., and van Oudenaarden A. (2005) Contributions of low molecule number and chromosomal positioning to stochastic gene expression. {\it Nature Genetics} 9, 937--944.

\mn
Beerenwinkel, N., Antal, T., Dingli, D., Traulsen, A., Kinzler, K.W., Velculescu, V.E., Vogelstein, B., and Nowak, M.A. (2007)
Genetic progression and the waiting time to cancer. {\it PLoS Computational Biology.} 3, paper e225

\mn
Beisel, C.J., Rokyta, D.R., Wichman, H.A., and Joyce, P. (2007) Testing the extreme value domain of attraction
for distributions of beneficial fitness effects. {\it Genetics.} 176, 2441--2449

\mn
Bodmer, W., and Tomlinson, I. (1995) Failure of programmed cell death and differentiation as causes of tumors: some simple mathematical models. {\it Proc Natl Acad Sci USA} 92, 11130--11134.

\mn
Coldman, A.J., and Murray, J.M. (2000) Optimal control for a stochastic model of cancer chemotherapy. {\it Mathematical Biosciences} 168, 187--200.

\mn
Cowperthwaite, M.C., Bull, J.J., and Meyers, L.A. (2005) Distributions of beneficial fitness effects in RNA.
{\it Gentics.} 170, 1449--1457

\mn
Durrett, R., and Mayberry, J. (2009) Traveling waves of selective sweeps.

\mn
Durrett, R., and Moseley, S. (2009) Evolution of resistance and progression to disease during clonal expansion
of cancer. {\it Theor. Pop. Biol.}, to appear

\mn
Elowitz, M.B. et al. (2002) Stochastic gene expression in a single cell. {\it Science} 297, 1183--1186.

\mn
Feinerman, O. et al. (2008) Variability and robustness in T cell activation from regulated heterogeneity in protein levels. {\it Science} 321, 1081.

\mn
Frank, S.A. (2007) \emph{Dynamics of Cancer: Incidence, Inheritance and Evolution.}
Princeton Series in Evolutionary Biology.

\mn
Gillespie, J.H. (1983) A simple stochastic gene substitution model. {\it Theor. Pop. Biol.} 23, 202--215

\mn
Gillespie, J.H. (1984) Molecular evolution over the mutational landscape. {\it Evolution.} 38, 1116--1129

\mn
Goldie, J.H., and Coldman, A.J. (1983) Quantitative model for multiple levels of drug resistance in clinical tumors. {\it Cancer Treatment Reports} 67, 923--931.

\mn
Goldie, J.H., and Coldman, A.J. (1984) The genetic origin of drug resistance in neoplasms: implications for systemic therapy. {\it Cancer Research} 44, 3643--3653.

\mn
Haeno, H., Iwasa, Y., and Michor, F. (2007) The evolution of two mutations during clonal expansion.
{\it Genetics.} 177, 2209--2221


\mn
Iwasa, Y., Michor, F., Komorova, N.L., and Nowak, M.A. (2005) Population genetics of tumor suppressor genes.
{\it J. Theor. Biol.} 233, 15--23

\mn
Iwasa, Y., Nowak, M.A., and Michor, F. (2006) Evolution of resistance during clonal expansion.
{\it Genetics.} 172, 2557--2566

\mn
Kassen, R., and Bataillon, T. (2006) Distribution of fitness effects among beneficial mutations before
selection in experimental populations of bacteria. {\it Nature Genetics.} 38, 484--488

\mn
Kaern, M. et al. (2005) Stochasticity in gene expression: from theories to phenotypes {\it Nature Reviews Genetics} 6, 451.

\mn
Knudson, A.D. (2001) Two genetic hits (more or less) to cancer. {\it Nature Reviews Cancer.} 1, 157--162

\mn
Komarova, N.L., and Wodarz, D. (2005) Drug resistance in cancer: principles of emergence and prevention. {\it Proc Natl Acad Sci USA} 102, 9714--9719.

\mn
Maley, C.C. et al. (2006) Genetic clonal diveresity predicts progression to esophageal adenocarcinoma.
{\it Nature Genetics.} 38, 468--473

\mn
Maley, C.C., and Forrest (2001) Exploring the relationship between neutral and selective mutations in cancer.
{\it Artif Life} 6, 325--345.

\mn
Michor, F., Iwasa, Y., and Nowak, M.A. (2004) Dynamics of cancer progression. {\it Nature Reviews Cancer}
4, 197--205

\mn
Michor, F., Nowak, M.A., and Iwasa, Y. (2006) Stochastic dynamics of metastasis formation. {\it J Theor Biol}
240, 521--530.

\mn
Michor, F., and Iwasa, Y. (2006) Dynamics of metastasis suppressor gene inactivation. {\it J Theor Biol}
241, 676--689.

\mn
Nowak M.A., Michor F, and Iwasa Y (2006) Genetic instability and clonal expansion. \emph{J Theor Biol} 241, 26--32.

\mn
Nowell P.C. (1976) The cloncal evolution of tumor cell populations. \emph{Science} 194, 23--28.

\mn
Orr, H.A. (2003) The distribution of fitness effects among beneficial mutations.
{\it Genetics.} 163, 1519--1526

\mn
Otto, S.P., and Jones, C.D. (2002) Detecting the undetected: Estimating the total number of loci
underlying a quantitative trait. {\it Genetics.} 156, 2093--2107

\mn
Rokyta, D.R., Beisel, C.J., Joyce, P., Ferris, M.T., Burch, C.L., and Wichman, H.A. (2008)
Beneficial fitness effects are not exponential in two viruses. {\it J. Mol. Evol.} 67, 368--376

\mn
Rozen, D.E., de Visser, J.A.G.M., and Gerrish, P.J. (2002) Fitness effects of fixed beneficial mutations
in microbial populations. {\it Curret Biology.} 12, 1040--1045

\mn
Sanju\'an, R., Moya, A., and Elena, S.F. (2004) The distribution of fitness effects caused by
single-nucleotide substitutions in an RNA virus. {\it Proc. Natl Acad. Sci., USA.} 101, 8396--8401

\mn
Shah, S.P., et al. (2009) Mutational evolution in a lobular breast tumour profiled at single
nucleotide resolution. {\it Nature.} 461, 809--813

\mn
Weissman, I. (1978) Estimation of parameters and large quantiles based on the $k$ largest observations.
{\it j. Amer. Stat. Assoc.} 73, 812--815

\mn
Wodarz, D., and Komarova, N.L. (2007) Can loss of apoptosis protect against cancer? {\it Trends Genet.} 23, 232--237.

\newpage

\end{document}